# Accurate calculation of oblate spheroidal wave functions

Arnie L. Van Buren

Cary, North Carolina

September 12, 2019

MSC -class: 3310

**Abstract**

Alternative expressions for calculating the oblate spheroidal radial functions of both kinds $R_{ml}^{(1)}(-ic,i\xi)$ and $R_{ml}^{(2)}(-ic,i\xi)$ are shown to provide accurate values over very large parameter ranges using double precision (64-bit) arithmetic, even where the traditional expressions fail. First is the expansion of the product of the radial function and the oblate angular function of the first kind $S_{ml}^{(1)}(-ic,\eta)$ in a series of products of spherical functions, with $\eta$ being a free parameter. Setting $\eta = 0$ leads to accurate values for $R_{ml}^{(2)}$ for $\xi \geq 0.01$ and non-small values of $l$ - $m$. Allowing it to vary with increasing $l$ leads to highly accurate values for $R_{ml}^{(1)}$ over all parameter ranges. Next is the calculation of $R_{ml}^{(2)}$ as an integral of the product of $S_{ml}^{(1)}$ and a spherical Neumann function kernel. This is useful for smaller values of $\xi$. Also used is the near equality of pairs of low-order eigenvalues when $c$ is large that leads to accurate values for $R_{ml}^{(2)}$ using neighboring accurate values for $R_{ml}^{(1)}$. A modified method is described that provides accurate values for the necessary expansion $d$ coefficients when $c$ is large and $l$ – $m$ is small and traditional methods fail. A resulting Fortran program oblfcn appears to provide useful values for the radial functions for virtually all desired values of $c$, $m$ and $\xi$ when using 64-bit arithmetic. Accuracies are almost always at least 8 decimal digits unless near a root or when $c$ is large, $\xi$ is near or equal to zero and $l$ - $m$ is small where the value of $R_{ml}^{(2)}$ or its first derivative is very small. Use of 128-bit arithmetic provides 15+ digits. Oblfcn is freely available at www.mathieuandspheroidalwavefunctions.com

## 1 Introduction

The scalar Helmholtz wave equation for steady waves, $(\nabla^2 + k^2)\Psi = 0$, where $k = 2\pi / \lambda$ and $\lambda$ is the wavelength, is separable in the oblate spheroidal coordinates $(\xi,\eta,\varphi)$, with $0 \leq \xi \leq \infty$, $-1 \leq \eta \leq 1$, and $0 \leq \varphi \leq 2\pi$. The factored solution is $\Psi_{ml}(\xi,\eta,\varphi) = R_{ml}(-ic,i\xi)\, S_{ml}(-ic,\eta)\, \Phi_m(\varphi)$, where $R_{ml}(-ic,i\xi)$ is the radial function, $S_{ml}(-ic,\eta)$ is the angular function, and $\Phi_m(\varphi)$ is the azimuthal function. Here $c = ka/2$, where $a$ is the interfocal distance of the elliptic cross section of the spheroid. The radial function of the first kind $R_{ml}^{(1)}(-ic,i\xi)$ and the radial function of the second kind $R_{ml}^{(2)}(-ic,i\xi)$ are the two independent solutions to the second order radial differential equation resulting from the separation of variables. These solutions are dependent on four parameters $(m,l,c,\xi)$ and an eigenvalue (separation constant) $A_{ml}(-ic)$. Similarly, $S_{ml}^{(1)}(-ic,\eta)$ and $S_{ml}^{(2)}(-ic,\eta)$ are the two independent solutions to the second order angular differential equation resulting from the separation of variables.

The addition of the imaginary number *i* to the notation for the oblate functions and eigenvalues reflects that the differential equations for oblate spheroidal functions can be obtained from the corresponding differential equations for prolate spheroidal functions by changing the size parameter *c* to *-ic* and the shape parameter $\xi$ to $i\xi$. Many of the expressions for calculating oblate functions can be obtained from prolate expressions using these same changes. In the following discussion we assume that the order *m* is either zero or a positive integer, with the degree *l* equal to *m, m+1, m+2, ….*

Oblate spheroidal functions are used in solving boundary value problems of radiation, scattering, and propagation of scalar and vector acoustic and electromagnetic waves in oblate spheroidal coordinates.

Existing computer programs for calculating the oblate radial functions including those developed at the Naval Research Laboratory [1-2] and those given by Zhang and Jin [3] utilize traditional expressions that fail to provide accurate results for some parameter ranges, especially at lower values of *l - m* when *c* becomes large and at large values of *l - m* for small to moderate values of $\xi$. Moreover, traditional methods fail to provide accurate values for the so-called *d* coefficients necessary for calculating both the oblate angular and radial functions at low values of *l – m* when *c* becomes large. Subtraction errors arising in evaluating the traditional expressions increase without bound as *c* increases. One approach to obtaining accurate function values in this case is to use arbitrary precision arithmetic, when available, and increase the precision to offset the increasing subtraction errors [4]. Of course, calculation times will be much greater than when using double precision arithmetic. The current paper describes several new expressions and procedures that extend the ability to calculate accurate radial function values over very wide parameter ranges using no more than double precision arithmetic.

A previous paper [5] presented an alternative expression for calculating prolate radial functions of the first kind $R_{ml}^{(1)}(c,\xi)$ that provides accurate values for all parameter ranges. The alternative expression follows directly from the product expansion provided by Meixner and Schäfke [6, p. 307] that represents the product of the radial and angular functions in a series of products of the corresponding spherical Bessel and associated Legendre functions. A second paper [7] presented alternative expressions that significantly increased the parameter ranges for which accurate values can be calculated for prolate radial functions of the second kind $R_{ml}^{(2)}(c,\xi)$. The first of these expressions utilizes the product expansion. The second expression was one given by Flammer [8, p. 53] that involves an integral of the product of $S_{ml}^{(1)}(c,\eta)$ and a spherical Neumann function kernel.

The present paper exploits the fact that the low-order eigenvalues for *c* large are paired. This near equality leads to accurate values for $R_{ml}^{(2)}$ and $R_{m,l+1}^{(2)}$ at low order using accurate values for neighboring $R_{m,l+1}^{(1)}$ and $R_{ml}^{(1)}$. The present paper presents a modified approach to calculating the oblate *d* coefficients that provides accurate values over virtually all parameter ranges, unlike the traditional method. It shows that the product expansion [6] leads to accurate calculation of $R_{ml}^{(1)}(-ic,i\xi)$ and improved calculation of $R_{ml}^{(2)}(-ic,i\xi)$. It shows that the integral expression given by Flammer provides accurate values for oblate radial functions in many parameter ranges where other expressions fail to do so. It also shows that an alternative Legendre function expansion for $R_{ml}^{(2)}(-ic,i\xi)$ due to Baber and Hasse [9] provides accurate values in some ranges where the traditional Legendre function expansion fails to do so. The paper then describes some of the features of a new Fortran computer program oblfcn that calculates the oblate spheroidal angular and radial functions using a combination of both the traditional and the alternative expressions. The paper concludes with a summary.

## 2 Angular functions of the first kind

The oblate angular function of the first kind $S_{ml}^{(1)}(-ic,\eta)$ is expressed [see for example ref. 8, p. 16] in terms of the corresponding associated Legendre functions of the first kind by

$$S_{ml}^{(1)}(-ic,\eta)=\sum_{n=0,1}^{\infty}{}' d_n(-ic\,|\,ml)\,P_{m+n}^{m}(\eta)\,, \tag{1}$$

where the prime sign on the summation indicates that $n=0,2,4,\ldots$ if $l-m$ is even or $n=1,3,5,\ldots$ if $l-m$ is odd. The following three term recursion formula relates successive expansion coefficients $d_{n-2}, d_n$, and $d_{n+2}$ for given values of *l, m,* and *c*:

$$\alpha_n d_{n+2}+(\beta_n-\lambda_{ml})\,d_n+\gamma_n d_{n-2}=0, \tag{2}$$

where

$$\alpha_n=-\frac{(2m+n+2)(2m+n+1)c^2}{(2n+2m+3)(2n+2m+5)},$$

$$\beta_n=\left[(n+m)(n+m+1)-\frac{2(n+m)(n+m+1)-2m^2-1}{(2n+2m+3)(2n+2m-1)}c^2\right],$$

$$\gamma_n=-\frac{n(n-1)c^2}{(2n+2m-3)(2n+2m-1)}. \tag{3}$$

Use of this formula to calculate the expansion coefficients requires a value for the separation constant or eigenvalue $\lambda_{ml}(-ic)$, which is chosen to ensure nontrivial convergent solutions for $S_{ml}^{(1)}(-ic,\eta)$.

For prolate functions the variational procedure developed by Bouwkamp [10] provides accurate values for each eigenvalue and its associated set of ratios of successive coefficients $d_{n+2}/d_n$. Starting with an estimate of the desired eigenvalue, one computes successive corrections to the estimate. Estimates for the first few eigenvalues are given by a power series expansion in $c^2$ when $c$ is less than about 8 and by an asymptotic expansion in $1/c$ for larger values of $c$. Estimates for higher order eigenvalues are given by extrapolation of lower order eigenvalues.

Usually only a few iterations are required for convergence to one of the prolate eigenvalues. If the estimate is close enough, the resulting eigenvalue will be the desired one $\lambda_{ml}$. Otherwise, convergence will be to one of the two neighboring eigenvalues of the same parity, $\lambda_{m,l-2}$ or $\lambda_{m,l+2}$, depending on whether the estimate is too low or too high . Use of the algorithm given in [6] assures that the procedure converges to the desired eigenvalue. Here an upper and lower bound for the interval containing the desired eigenvalue are established based on previous eigenvalues. The interval is chosen to ensure it does not contain either of the neighboring eigenvalues of the same parity. If the Bouwkamp procedure converges to a value within the interval, it is the desired eigenvalue. If it converges to the eigenvalue below the interval $\lambda_{m,l-2}$, then the estimate is taken as a new lower bound. If it converges to the eigenvalue above the interval $\lambda_{m,l-2}$, the estimate is taken as a new upper bound. A new estimate is set equal to the midpoint of the new bounds and the procedure is repeated. The estimate continues to become more accurate and only a few iterations are required until the estimate is accurate enough for full convergence of the Bouwkamp procedure, even for those rare cases where an estimate accurate to 5 or 6 digits is insufficient.

For oblate functions the situation is more complicated. When *c* is less than about 40, the Bouwkamp procedure converges to the desired eigenvalue for all values of *m* and *l* when the estimate is accurate to at least one part in $10^5$. For larger values of *c*, the Bouwkamp procedure works well when $l-m$ is somewhat larger than the so-called breakpoint $nb$, equal to $2c/\pi$. For lower values of $l-m$, however, the eigenvalue estimate must be increasingly accurate as *c* becomes larger. If the estimate is not accurate enough, the Bouwkamp procedure will still converge to an eigenvalue of the same parity, but one corresponding to a much higher value of *l* - *m*. For example, the required accuracy for *c* = 100 is about one part in $10^{12}$; for *c* = 200 about one part in $10^{24}$; for *c* = 500 even 32 digits of accuracy is insufficient. Thus an alternative matrix procedure is used to obtain accurate lower order eigenvalues.

By successively choosing *n* = 0, 1, 2, 3,... , in (2), one can obtain an infinite set of simultaneous equations for the coefficients $d_n$. These equations can be written in matrix form as

$$\{B\}\{d\} = \lambda_{ml}\{d\}, \tag{4}$$

where $\{B\}$ is an infinite square tridiagonal matrix depending on *m* and *c*, $\{d\}$ is a vector representation of the $d_n$ coefficients, and $\lambda_{ml}$ is the eigenvalue for *m* and *l*. The desired oblate eigenvalues are then the set of eigenvalues of $\{B\}$.

King and Hanish [11] show that the matrix becomes symmetrical when $d_r$ is replaced with $d_r = \sqrt{\frac{(2r+2m+1)(r!)}{2(r+2m)!}}\, d_r$ , *r* = *n-2, n,* and *n+2*. This is equivalent to using associated Legendre functions with unit normalization in (1). It is much easier and faster to compute eigenvalues of a symmetric tridiagonal matrix. Furthermore, the $d_n$ coefficients with even subscript are only involved in spheroidal functions with even $l-m$, while those with odd subscript are involved when $l-m$ is odd. This allows the matrix $\{B\}$ to be decomposed into an even matrix $\{B^e\}$ using *n* = 0, 2, 4,..., whose eigenvalues are for $l-m$ even and an odd matrix $\{B^o\}$ using *n* = 1, 3, 5,..., whose eigenvalues are for $l-m$ odd. It is convenient to divide all of the matrix elements by $c^2$ and truncate both matrices to either order $4n_b/3$ or order 67 ,whichever is larger. Use of a standard tridiagonal matrix routine results in an odd and an even set of eigenvalues. Ordering each set of eigenvalues in increasing numerical value and interlacing the two sets results in accurate oblate eigenvalues $\lambda_{ml}(-ic)$ for $l-m$ =0, 1, 2, ..., $n_b$.

Although the matrix results are nearly fully accurate for values of $l-m$ up to $n_b$, the Bouwkamp procedure is attempted for all values of $l-m$ regardless of the value for *c*. This often provides eigenvalues for lower values of $l-m$ that are slightly more accurate than the matrix results. The matrix values are used as starting values for $l-m$ up to $4n_b/3$ and estimates using extrapolation from previous eigenvalues are used for higher values of $l-m$. For lower values of $l-m$,when the Bouwkamp procedure fails to converge to an eigenvalue close to the matrix value, the matrix value is taken as the eigenvalue.

An accurate eigenvalue allows one to now compute the $d_n$ coefficients using the recursion formula (2). Dividing each term in (2) by $d_n$ results in an expression relating the ratio $N_{n+2} = d_{n+2}/d_n$ to the ratio $N_n = d_n/d_{n-2}$. Traditionally this expression is used in the forward direction to calculate ratios up to $N_{l-m}$ starting with the first ratio $N_2 = (\beta_0 - \lambda_{ml})/\alpha_0$ for *n* even or $N_3 = (\beta_1 - \lambda_{ml})/\alpha_1$ for *n* odd. Ratios for *n* above $l-m$ are calculated backward from a suitably high value of *n* where the ratio $N_{n+2}$ is set equal to zero. As the ratios are calculated backward, they become progressively more accurate until they are essentially fully accurate. Only the backward recursion is used when $l-m=0$ or 1. See Flammer (7) for a discussion of this. This procedure always works well for the prolate case. Here the forward recursion provides accurate values up the ratio $N_{l-m}$ and the backward recursion provide accurate values down to the ratio $N_{l-m+2}$. This is not always true in the oblate case.

Here, for high values of *c* and lower values of $l-m$, the values of *n* for which both the forward and the backward recursion maintain accuracy is no longer $l-m$ but instead is a different value $n^*$ . When $n^*$ is smaller than $l-m$ this results is a loss of accuracy in ratios calculated using the forward recursion starting with $N_{n^*+2}$. The loss of accuracy increases until $N_{l-m}$ is reached. The backward recursion must be used here to obtain accurate values for ratios from $N_{l-m}$ down to $N_{n^*}$ Similarly, when $n^*$ is larger than $l-m$, there is a loss of accuracy in ratios calculated using the backward recursion starting with $N_{n^*}$. The loss of accuracy increases until the ratio $N_{l-m+2}$ is reached. The forward recursion must be used here to obtain the ratios $N_{l-m+2}$ to $N_{n^*}$.

The following procedure provides accurate ratios for the oblate case. One first uses the backward recursion down to the either $N_2$ for *n* even or $N_3$ for *n* odd. Then one uses the forward recursion from either $N_2$ or $N_3$ until the forward and backward ratio values match to ndec - 1 decimal digits, where ndec is the number of decimal digits available in floating point numbers. If there is no match to ndec - 1 digits, the forward recursion is continued until the match starts decreasing significantly. The best match is then selected. It is rarely less than ndec - 2 digits.

If one uses the traditional procedure instead, many of the resulting coefficient ratios will be less than fully accurate when *c* is large and $l-m$ is smaller than $n_b$. This will then reduce the accuracy of the resulting oblate function values. For a given large value of *c*, the match at the traditional value of *n* equal to $l-m$ tends to be high for m = 0 but decreases as *c* increases until it reaches a minimum at a value of m approximately equal to 0.3*c*. The minimum always occurs at or near $l = m$. Figure 1 shows the minimum number of decimal digits in the match at $l = m$ as a function of c. For $c$ greater than 300, the minimum is zero digits.

The coefficients can be normalized by requiring that $S_{ml}^{(1)}(-ic,\eta)$ has the same normalization factor as $P_l^m(\eta)$ [5], resulting in the following relation:

$$\sum_{n=0,1}^{\infty}{}' \frac{2(n+2m)!}{[2(n+m)+1]n!}[d_n(-ic\,|\,ml)]^2 = \frac{2(l+m)!}{(2l+1)(l-m)!}\,. \qquad (5)$$

Use of this has the practical advantage of eliminating the need to compute the normalization factor which is often found in problems involving expansions in spheroidal angular functions. A better choice, however, is to set the right-hand side of (5) equal to unity. This results in the angular functions having unit norm. It has the advantage of limiting the magnitude of the angular functions to moderate values. For other normalization schemes such as (5) the angular functions become increasing large in magnitude as $m$ increases and eventually overflow even in 128-bit arithmetic with a maximum exponent over 4000. The program oblfcn offers either unit normalization or the same normalization as the corresponding Legendre functions.

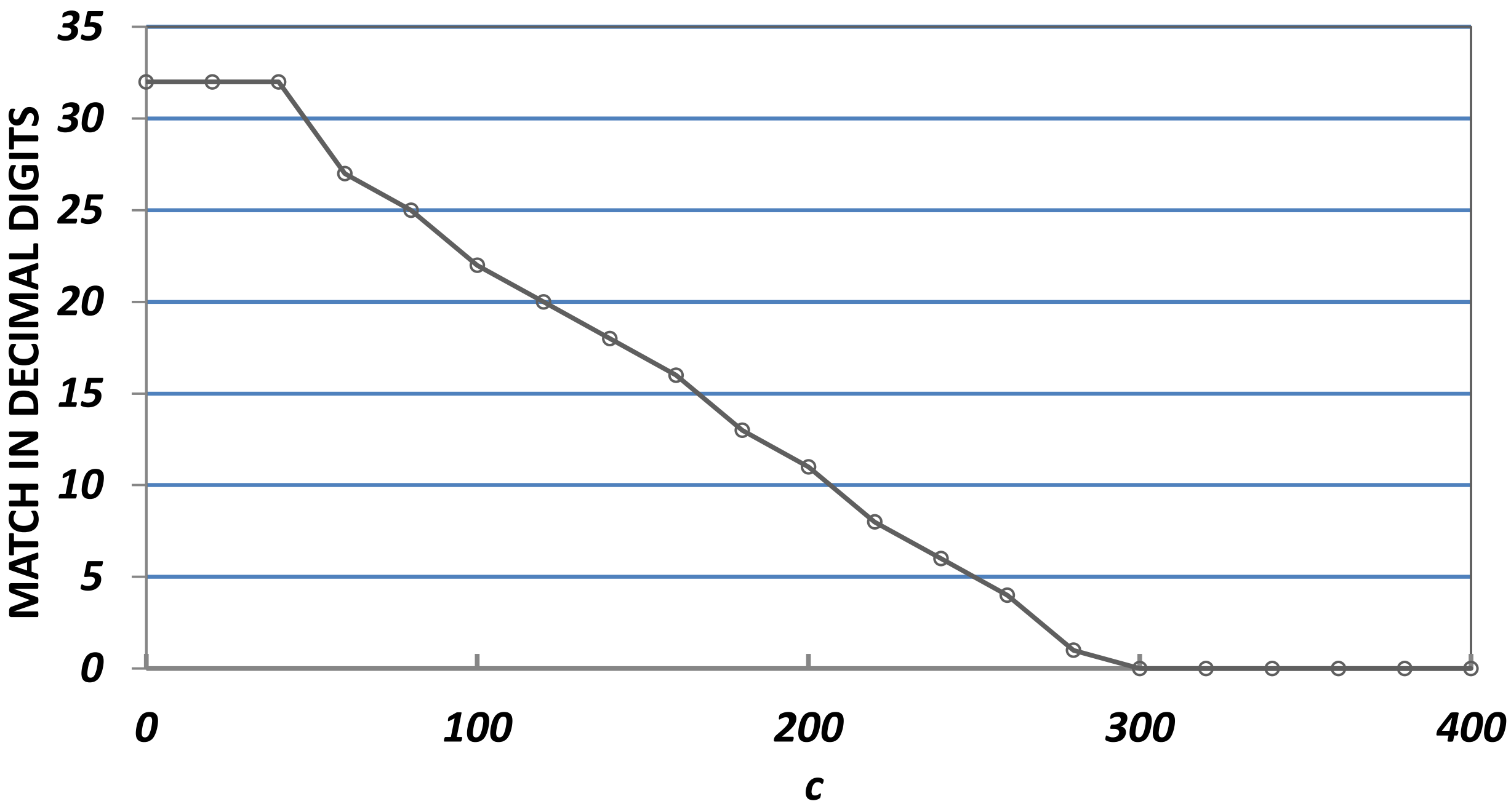

Fig 1: Minimum number of matching decimal digits for the value of $d_{l-m+2}/d_{l-m}$ computed from forward recursion to its value computed from backward recursion plotted versus $c$.

The normalization sum in (5) is numerically robust with no subtraction errors occurring in its computation. This contrasts with the corresponding normalization sum involved in requiring the angular functions to match the corresponding associated Legendre functions at $\eta = 0$ [see e.g., Flammer [8], p. 21]. Here subtraction errors can occur in calculating the angular function at $\eta = 0$, especially when $c$ is large and $l-m$ is less than the breakpoint $nb$. For fixed $c$ the error is maximum at $m = 0$ and for fixed $l-m$ decreases to zero as $m$ increases. For fixed $c$ and $m$, it is maximum at $l = m$ and decreases to zero as the breakpoint is approached. Figure 2 shows the subtraction error in decimal digits at $l = m$ plotted versus $c$ for several values of $m$.

Calculation of the angular functions using (1) can suffer subtraction errors at large values of $c$ for values of $\eta$ other than 0. For a given value of $m$ and for $l - m$ less than $n_b$, the error decreases to zero as $\eta$ increases from 0. For given values of $m$ and $\eta$, the error decreases to zero as $l-m$ increases. Figure 3 shows several examples of the subtraction error in decimal digits occurring in calculating the angular function $S_{ml}^{(1)}(-ic,\eta)$, plotted as a function of *c*. There is no subtraction error involved in calculating the angular functions when $\eta$ is equal to 0, regardless of the value for *c*.

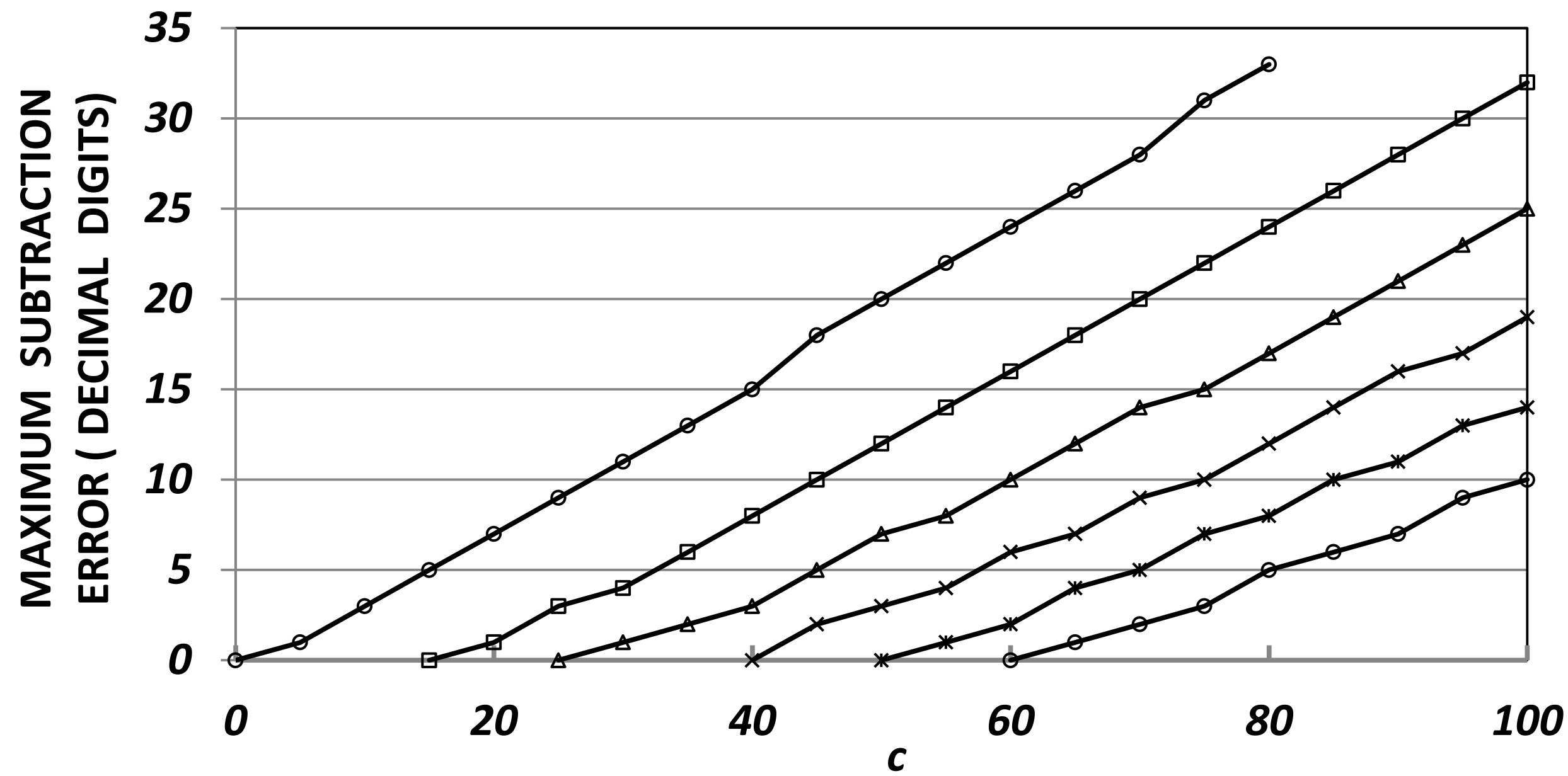

Fig 2: Subtraction error in decimal digits involved in calculating $S^{(1)}_{m,m}(-ic,0)$ plotted versus $c$ for selected values of m: o (0), □ (10), Δ (30), × (40), * (50), ◊ (60).

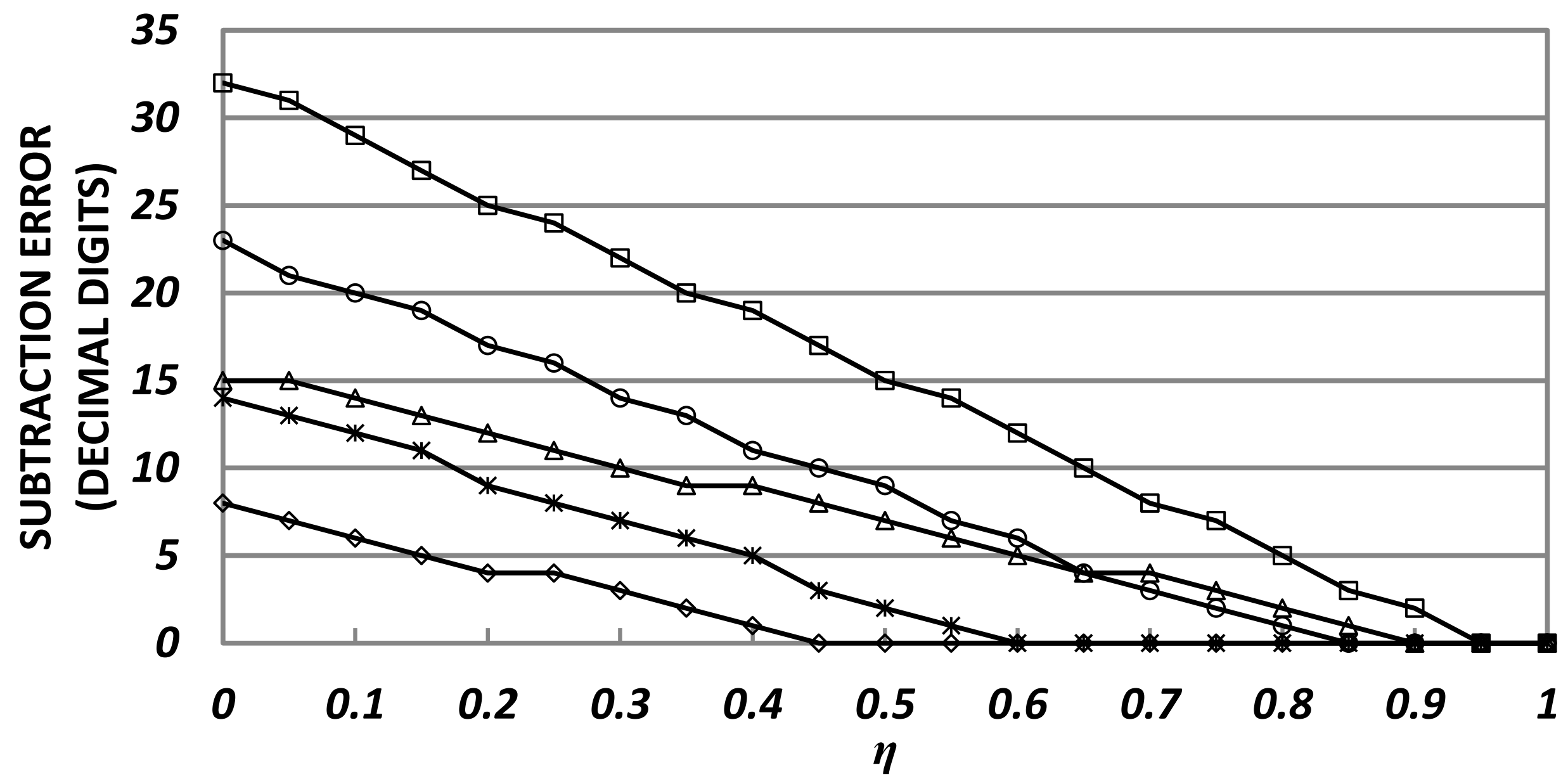

Fig. 3: Subtraction error when calculating angular functions of the first kind for selected parameters ($c, m ,l$): Δ (40, 0, 0); □ (80, 0, 0); * (80, 0, 20); o (80, 10, 10); ◊ (80, 10, 30).

## 3 Expansion of the product of the radial and angular functions

The expansion of the product of $R_{ml}^{(j)}(-ic,i\xi)$ and $S_{ml}^{(1)}(-ic,\eta)$ in terms of the corresponding spherical functions is given by:

$$R_{ml}^{(j)}(-ic,i\xi)\,S_{ml}^{(1)}(-ic,\eta)=\sum_{n=0,1}^{\infty}{}'\,i^{n+m-l}\,d_n(-ic\,|\,ml)\,\psi_{n+m}^{(j)}(kr)\,P_{n+m}^{m}(\cos\theta)\,, \tag{6}$$

where $j$ = 1 or 2. $\psi_{m+n}^{(1)}(kr)$ is the spherical Bessel function $j_{m+n}(kr)$ and $\psi_{m+n}^{(2)}(kr)$ is the spherical Neumann function $y_{n+m}(kr)$. This is a special case of the more general expansion given by Meixner and Schäfke [6, p. 307]. Using the relationship between the spherical coordinates *r* and $\theta$ and spheroidal coordinates (about the same origin and with $\eta$ = 1 coincident with $\theta$ = 0) we obtain $kr=c(\xi^2-\eta^2+1)^{1/2}$ and $\cos\theta=\eta\xi/(\xi^2-\eta^2+1)^{1/2}$. Substituting for $S_{ml}^{(1)}(-ic,\eta)$ from (1) and solving for $R_{ml}^{(j)}(-ic,i\xi)$ produces

$$R_{ml}^{(j)}(-ic,i\xi)=\frac{\displaystyle\sum_{n=0,1}^{\infty}{}'\,i^{n+m-l}\,d_n(-ic\,|\,ml)\,\psi_{n+m}^{(j)}[c(\xi^2-\eta^2+1)^{1/2}]\,P_{n+m}^{m}[\eta\xi/(\xi^2-\eta^2+1)^{1/2}]}{\displaystyle\sum_{n=0,1}^{\infty}{}'\,d_n(-ic\,|\,ml)\,P_{n+m}^{m}(\eta)}\,. \tag{7}$$

The significance of this general expression is that it allows us to choose the value for $\eta$ that provides the maximum accuracy for calculated values of $R_{ml}^{(j)}(-ic,i\xi)$. For many parameter ranges it is desirable to allow $\eta$ to vary as the value of the index *l* increases from *m* to higher values. Application of this so-called variable $\eta$ method to the calculation of the radial functions of both the first and second kinds will be described below.

## 4 Traditional Bessel function expressions

Consider the case when $\eta$ = 1. The argument of $P_{n+m}^{m}$ in both the numerator and the denominator approaches unity as $\eta$ approaches unity. Although $P_{n+m}^{m}$ approaches zero in this case for $m\neq 0$, the limit of the rhs of (7) exists and we obtain:

$$R_{ml}^{(j)}(-ic,i\xi)=\left(\frac{\xi^2+1}{\xi^2}\right)^{m/2}\frac{\displaystyle\sum_{n=0,1}^{\infty}{}'\,i^{n+m-l}\,d_n(-ic\,|\,ml)\,\psi_{m+n}^{(j)}(c\xi)\,\frac{(n+2m)!}{n!}}{\displaystyle\sum_{n=0,1}^{\infty}{}'\,d_n(-ic\,|\,ml)\,\frac{(n+2m)!}{n!}}\,. \tag{8}$$

Flammer [8, p. 32] derives (8) using integral representations of the spheroidal wave functions. The corresponding expression for the first derivative of $R_{ml}^{(j)}$ with respect to $\xi$ is obtained by taking the first derivative of the rhs of (8). Equation (8) is the expression commonly used to calculate numerical values for both $R_{ml}^{(1)}$ and $R_{ml}^{(2)}$. The advantage of these expressions for the oblate case is that the denominator sum is robust with no subtraction errors.

## 5 Alternative Bessel function expressions

Reference [5] also provides expressions for the prolate radial functions of the first kind obtained by choosing $\eta = 0$ in the prolate version of (7). Converting these expressions to oblate form and extending them to include radial functions of the second kind results in the following alternative expressions:

$$R_{ml}^{(j)}(-ic, i\xi) = \frac{\sum_{n=0,1}^{\infty}{}' i^{n+m-l}\, d_n(-ic\,|\,ml)\, \psi_{m+n}^{(j)}(c\sqrt{\xi^2+1})\, P_{n+m}^m(0)}{\sum_{n=0,1}^{\infty}{}' d_n(-ic\,|\,ml)\, P_{n+m}^m(0)}, l-m\text{, even} \quad (9)$$

and

$$R_{ml}^{(j)}(-ic, i\xi) = \left(\frac{\xi}{\sqrt{\xi^2+1}}\right) \frac{\sum_{n=0,1}^{\infty}{}' i^{n+m-l}\, d_n(-ic\,|\,ml)\, \psi_{m+n}^{(j)}(c\sqrt{\xi^2+1}) \left[\frac{dP_{n+m}^m(\eta)}{d\eta}\right]_{\eta=0}}{\sum_{n=0,1}^{\infty}{}' d_n(-ic\,|\,ml) \left[\frac{dP_{n+m}^m(\eta)}{d\eta}\right]_{\eta=0}}, l-m\text{ odd.} \quad (10)$$

Reference [5] shows that the corresponding expressions for the prolate case are numerically robust when calculating the prolate radial function of the first kind. They suffer no subtraction errors and provide accurate values over all parameter ranges. For the oblate case the $\eta = 0$ expressions will prove to be useful in the calculation of the radial functions of the second kind.

## 6 Calculation of the oblate radial functions of the first kind

The traditional expressions (7) and (8) provide accurate values for the oblate radial functions of the first kind unless *c* is not small and *l* - *m* is less than the break point $n_b$. In this case there are subtraction errors that are maximum at *l* - *m* = 0 and that decrease monotonically to zero with increasing l - m until the breakpoint is reached. Subtraction error is defined to be the number of accurate decimal digits that are lost in calculating the sum of the series. This loss of accuracy occurs when the sum of all of the positive terms in the series is nearly equal to the sum of all of the negative terms. The subtraction error is then equal to the number of leading decimal digits that are the same in the positive and negative sums. Figure 4 shows examples of the subtraction error arising when calculating $R_{ml}^{(1)}$ and its first derivative using the traditional Bessel function expressions (7) and (8). The subtraction errors are zero or nearly so for *m* = 0, increase with increasing *m* until *m* is large and then decrease with increasing *m*. They also tend to be largest when $\xi$ is near unity.

These subtraction errors can be avoided by using the variable $\eta$ method to calculate $R_{ml}^{(1)}$ and its first derivative. If the subtraction errors for $l = m$ are not near zero using the traditional expressions, then the value for $\eta$ is reduced in steps from unity until the errors are as small as possible. Usually both subtraction errors will now be near zero. However, if either $R_{ml}^{(1)}$ or its first derivative is near a root, its value will be somewhat smaller in magnitude than expected and its calculation will involve an unavoidable subtraction error whose size depends on how close the root is.

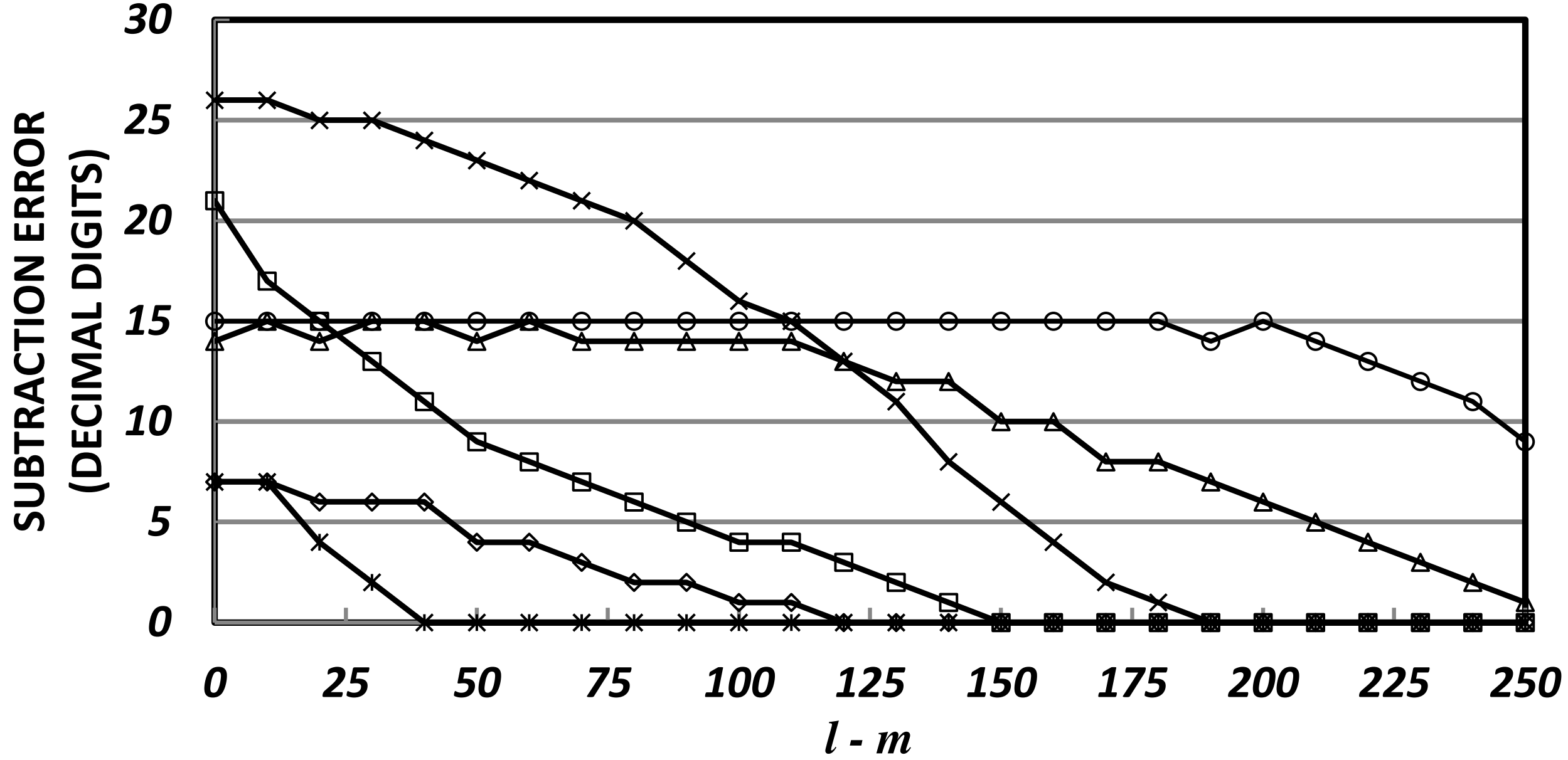

Fig. 4: Subtraction error when calculating $R_{ml}^{(1)}(-ic, i\xi)$ and its first derivative using the traditional Bessel functions expressions for selected parameters ($\xi$, $c$, $m$): ◊ (0.1, 300, 50); □ (0.5, 300, 100); Δ (1.0, 300, 100); × (1.0, 300, 200); * (1.0, 100, 100); o (1.5, 300, 200).

Using equal steps in $\theta = \arccos(\eta)$ works well, with a step size $\Delta\theta$ of about 0.05 radian. Here $\theta$ is incremented from zero, the value for the traditional expressions. Only a few steps are usually required. The value of $\theta$ that worked for $l = m$ is then used for progressively higher values of $l$ until the subtraction error increases. Then the value of $\theta$ is incrementally decreased until the subtraction error is near zero again. Usually only one step is required. The process is continued either until the maximum desired value of $l$ is reached or until $\theta$ has reached zero. Once $\theta$ reaches zero, the traditional expressions work well for higher values of $l$. Figure 5 shows the number of steps in $\theta$ that provide nearly zero subtraction error for several of the examples given in Fig. 4. For all but the case where $\xi = 1.5$, the traditional expressions also work well for $l$ - $m$ greater than 250. The variable $\eta$ algorithm had not yet decreased $\theta$ to zero when $l$ - $m$ reached 250, since the $\theta$ values shown in Fig. 5 were still providing subtraction errors near zero.

## 7 Use of low-order pairing of eigenvalues at higher values of c

Calculation of the radial functions of the second kind can be difficult. This is especially true at low orders when $c$ is large and $\xi$ is somewhat less than unity. Traditional expansions often fail here. Fortunately one can take advantage of the low-order pairing of eigenvalues. For non-small values of $c$, the eigenvalue $A_{mm}$ is nearly equal to $A_{m,m+1}$; $A_{m,m+2}$ is nearly equal to $A_{m,m+3}$; .... The agreement decreases with increasing $l$ and disappears as $l$ - $m$ approaches 2c/π. The differential equations for the corresponding low-order radial functions are then nearly identical as are their solutions. The following approximations are found to apply here:

$$\begin{aligned} R_{m,l}^{(2)}(-ic,i\xi) &\simeq R_{m,l+1}^{(1)}(-ic,i\xi),\, l\ even, \\ &\simeq -R_{m,l-1}^{(1)}(-ic,i\xi),\, l\ odd. \end{aligned} \tag{11}$$

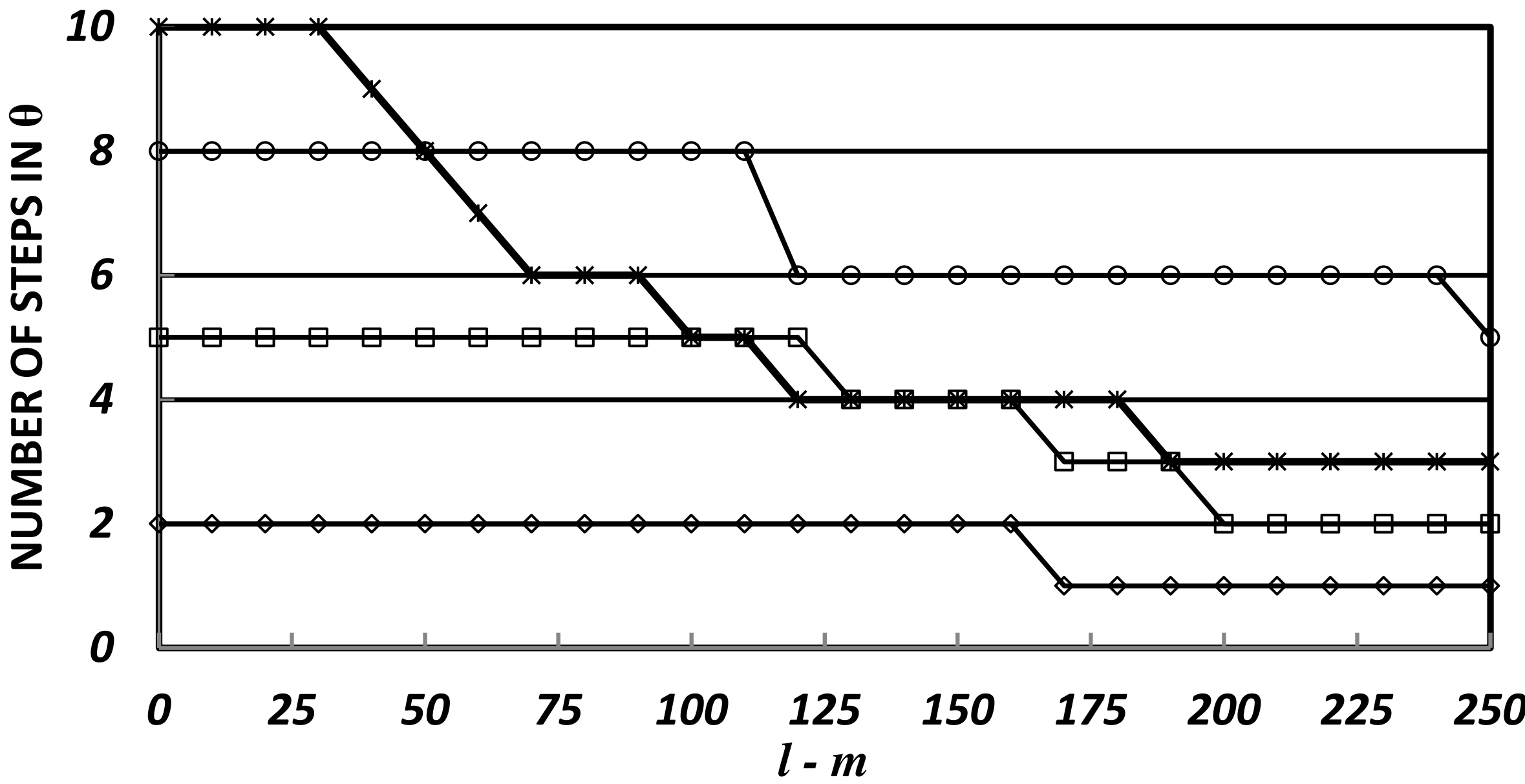

Fig. 5: Number of steps in $\theta$ required to calculate $R_{ml}^{(1)}(-ic,i\xi)$ and its first derivative with virtually no subtraction error for selected parameters ($\xi$, $c$, $m$): ◊ (0.1, 300, 50); □ (0.5, 300, 100); * (1.0, 100, 100); o (1.5, 300, 200).

This is equivalent to setting $R_{m,m+1}^{(3)} = R_{m,m+1}^{(1)} + iR_{m,m+1}^{(2)} \simeq -i\left[R_{m,m}^{(1)} + iR_{m,m}^{(2)}\right]$ and so forth for higher values of $l$ with paired eigenvalues. Flammer [8, p.68] obtained these approximations for the special case of $\xi = 0$ based on an asymptotic representation of the radial function $R_{m,l}^{(3)}$ in a series of Laguerre functions. Observation of numerical results show that (11) is equally valid for other values of $\xi$.

The agreement between neighboring eigenvalues for $c$ = 400 and $m$ = 0. 50 and 100 is shown in Fig. 6. Calculations were performed in 128 bit arithmetic with 33 decimal digits available. Fig. 6 also shows the agreement between the corresponding radial functions in accordance with (11) for the case where $\xi$ = 1.5. Both $R_{m,l}^{(1)}$ and $R_{m,l}^{(2)}$ plus their first derivatives are compared and the minimum of the four results is plotted. Here $R_{ml}^{(2)}$ and its first derivative have been calculated using Neumann function expansions that provide nearly full accuracy for large c and $\xi$ greater than about unity The agreement at lower values of *l* - *m* is limited by the precision available. Also the function values are expected to be less accurate than their eigenvalues. If much higher precision were used, the agreement would certainly be corresponding higher. A conservative estimate of the accuracy of values for $R_{ml}^{(2)}$ and its first derivative obtained using (11) is given by subtracting 2 from the corresponding eigenvalue agreement.

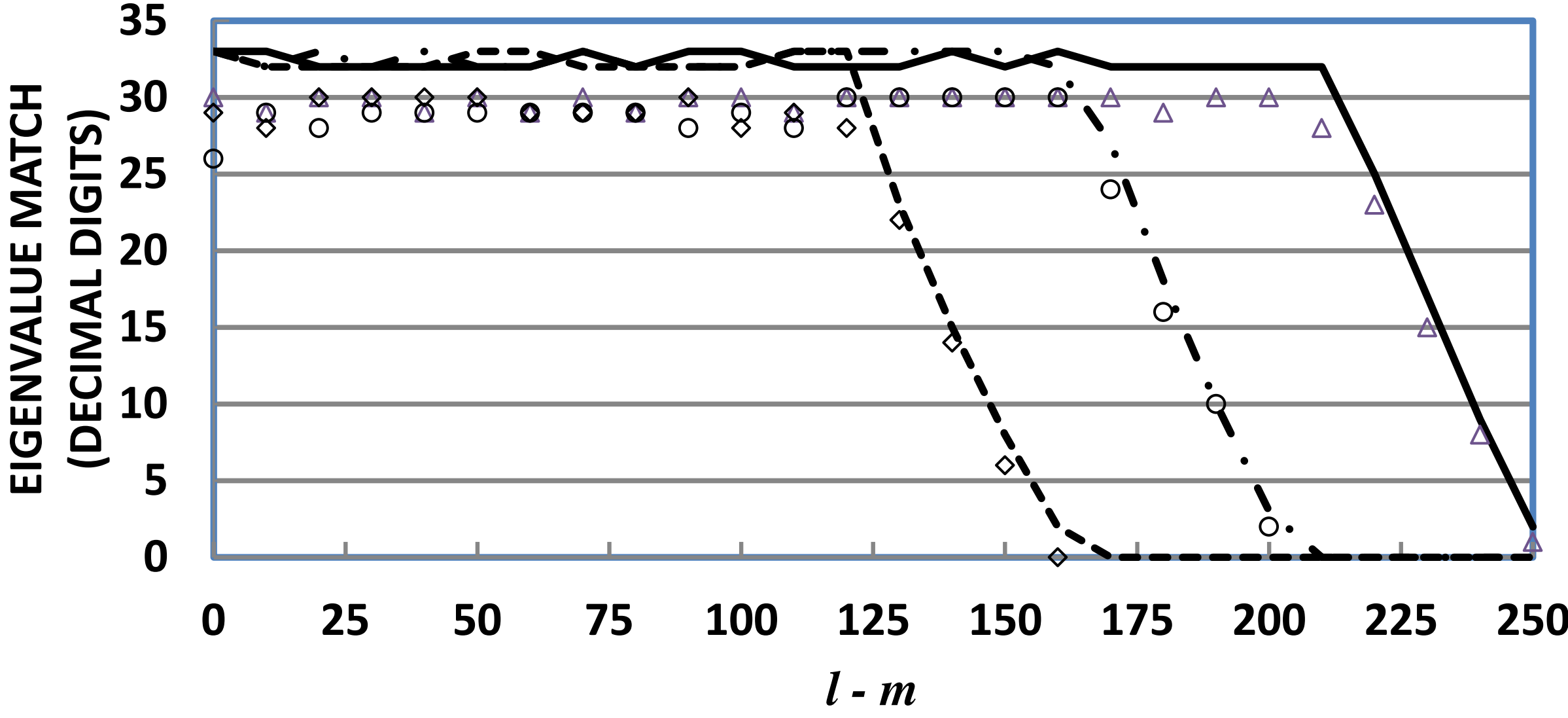

Fig. 6: Agreement in decimal digits between neighboring low-order eigenvalues for $c$ = 400 and $m$ = 0 (solid line), $m$ = 50 (dot dash line) and $m$ = 100 (dash line). Minimum agreement between corresponding radial function values in accordance with (11) for $\xi$ = 1.5 for $m$ = 0 (Δ), $m$ = 50 (o), and $m$ = 100 (◊).

## 8 Traditional Legendre function expression

A variety of methods are required to obtain accurate values for $R_{ml}^{(2)}$ over a wide range of parameters. Considered next is use of the traditional Legendre function expression {see, e.g., Flammer [8, p. 33] or Zhang and Jin [3, p. 567]}. Reference [2] shows that this expression can be obtained using the product expansion. Here the angular function of the first kind is replaced with that of the second kind and the associated Legendre function of the first kind $P_{n+m}^{m}$ is replaced with the associated Legendre function of the second kind $Q_{n+m}^{m}$. The resulting expression is:

$$R_{ml}^{(2)}(-ic,i\xi) = \frac{1}{\kappa_{ml}^{(2)}(-ic)} \sum_{n=-\infty}^{\infty}{}' \, d_n(-ic \mid ml) Q_{n+m}^{m}(i\xi). \tag{12}$$

Of course, equation (12) could have been obtained more directly from the fact that $R_{ml}^{(2)}$ and $S_{ml}^{(2)}$ are proportional to each other.

Although $Q^m_{n+m}(i\xi)$ becomes infinitely large when $n$ is less than -2$m$, its product with $d_n$ is finite and proportional to $P^m_{-n-m-1}(i\xi)$. The rhs of (12) then divides into two series, one over $n$ from -2$m$ (or -2$m$ + 1 if $l$ - $m$ is odd) to ∞ involving $Q^m_{n+m}(i\xi)$ and one over $n$ from 2$m$ + 2 (or 2$m$ + 1 if $l$ - $m$ is odd) to ∞ involving $P^m_{n-m-1}(i\xi)$. The result is the traditional associated Legendre function expression used to evaluate $R^{(2)}_{ml}(-ic,i\xi)$ when $\xi$ is small. The $d_n$ coefficients with negative subscripts required in the $Q^m_{n+m}(i\xi)$ series and the special $d_{\rho/n}$ coefficients required for the $P^m_{n-m-1}(i\xi)$ series can be computed from either $d_0$ or $d_1$, depending on whether $l$ - $m$ is even or odd. The joining factor $\kappa^{(2)}_{ml}(-ic)$ contains the same series that is given in the denominator of the traditional Bessel function expression (8). Thus the denominator of the associated Legendre function expression (12) suffers a loss in accuracy for low values of $l$ - $m$ similar to that shown in Fig. 2 due to subtraction errors that increase without bound as $c$ increases. The Legendre functions sums in (12) suffer even greater subtraction error than the denominator does, especially for low values of $m$. Figure 7 shows examples of the subtraction error in the Legendre function sums for selected parameter sets. The subtraction error for $\xi$ less than 0.5 is nearly the same as that shown for ξ = 0.5. The error decreases with increasing $l$ - $m$, reaching zero near the breakpoint for low $m$ and somewhat beyond for larger $m$. For a fixed value of $l$ -$m$, the error tends to increase as $\xi$ becomes larger than 0.5, increases with increasing $c$ and decreases with increasing $m$. The traditional Legendre function expression is used when $\xi$ is less than or equal to 0.99. Other methods will be used to compute $R^{(2)}_{ml}$ and its first derivative when $\xi$ is larger or when the subtraction error in the traditional Legendre function expression is too large to allow sufficiently accurate values.

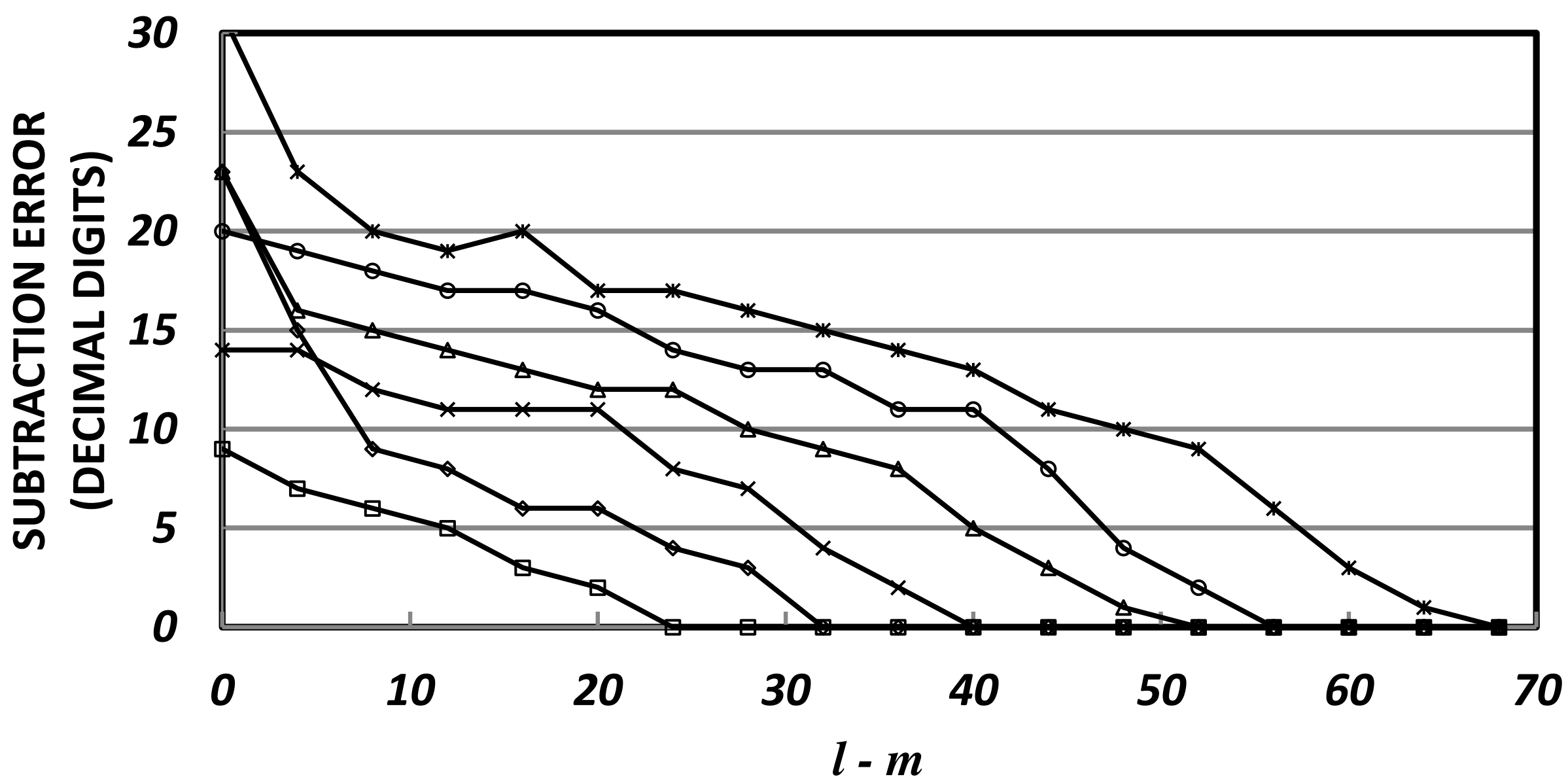

Fig. 7: Subtraction error in decimal digits occurring in the traditional Legendre function series for selected parameters ($\xi$, $c$, $m$): ◊ (0.5, 30, 0); □ (0.5, 30, 10); Δ (0.99, 30, 0); × (0.99, 30, 10); * (0.99, 40, 0); o (0.99, 40, 10).

The accuracy of the set of radial function values for $R_{ml}^{(1)}$, $R_{ml}^{(2)}$ and their first derivatives for a given value of $m$ and $l$ can be estimated using the Wronskian relationship:

$$R_{ml}^{(1)}\frac{dR_{ml}^{(2)}}{d\xi}-R_{ml}^{(2)}\frac{dR_{ml}^{(1)}}{d\xi}=\frac{1}{c(\xi^2+1)}. \tag{13}$$

An integer estimate of accuracy is given by the number of leading digits of agreement between the calculated Wronskian on the left-hand side of (13) and the theoretical Wronskian given by the right-hand side of (13).

The traditional Bessel function expressions together with the variable $\eta$ method provide accurate values for $R_{ml}^{(1)}$ and its first derivative. Thus the number of decimal digits of agreement between the theoretical and computed Wronskian is a measure of the accuracy of $R_{ml}^{(2)}$ and its first derivative.

## 9 Alternative Legendre function expression

Baber and Hasse [9] provided the following expression for the oblate radial functions of the third kind $R_{ml}^{(3)} = R_{ml}^{(1)} + iR_{ml}^{(2)}$ in terms of the functions $Q_{n+m}^{m}(i\xi)$:

$$R_{ml}^{(3)}(-ic,i\xi)=\frac{e^{ic\xi}i^{2m-l}}{m!c}\sum_{n=-l}^{\infty}\frac{A_n^{ml}}{A_{-m}^{ml}}Q_{m+n}^{m}(i\xi). \tag{14}$$

The coefficients $A_n^{ml}$ satisfy the following recursion relation:

$$\begin{aligned}&\frac{2c(n+m+1)(n+2m+1)}{(2n+2m+3)}A_{n+1}^{ml}-[(n+m)(n+m+1)-\lambda_{ml}-c^2]A_n^{ml}\\&-\frac{2cn(n+m)}{(2n+2m-1)}A_{n-1}^{ml}=0,\end{aligned} \tag{15}$$

with the asymptotic condition $A_{n+1}^{ml}/A_n^{ml}\xrightarrow[n\to\infty]{}c/n$. Flammer [8, p. 40] and ref. [2] provide a discussion of this expression. The radial functions of the second kind are then given by the imaginary part of the right hand side of (14). This expression does not suffer subtraction errors for the lowest order functions when both $c$ and $m$ are small to moderate in size and $\xi$ is not large. It can provide accurate results at lower values of $l$ - $m$ when the traditional Legendre function expression fails to do so. Examples of this are when $c$ = 60, $m$ = 0 and $\xi$ is either 0.1 or 0.5. Figure 8 shows the Wronskian accuracy in decimal digits for values of $R_{ml}^{(2)}$ and its first derivative obtained using both the traditional Legendre expression and equation (14). The calculations were performed in 128-bit arithmetic with 33 decimal digits of precision. The two methods complement each other well, with one or both providing at least 16 digits of accuracy over the full range of $l$ - $m$. Similar results can be obtained with $m$ as large as 40. Backward recursion from a suitably large index, where the ratio is set equal to 0, down to $n$ = 0 provides ratios of successive expansion coefficients $A_{n+1}^{ml}/A_n^{ml}$ for n ≥ 0. When m > 0, the ratio $A_{-m+1}^{ml}/A_{-m}^{ml}$ is obtained from (15) by setting $n$ = $-m$ where the third term is identically 0. Forward recursion to $n$ = -1 then produces the remaining ratios. Multiplication by successive ratios starting with $A_{-m+1}^{ml}/A_{-m}^{ml}$ times $A_{-m+2}^{ml}/A_{-m+1}^{ml}$ gives the coefficients used in (14).

## 10 Integral expressions for calculating radial functions

Reference [6] shows that the integral expressions given by Flammer [8, pp.53-54] are useful for calculating the prolate radial functions of the second kind. Converted to oblate form, these expressions become:

$$R_{ml}^{(2)}(-ic,i\xi)=\frac{(-1)^{(l-m)/2}(2m+1)}{2^{m+1}m!d_0(-ic\,|\,ml)}\times$$
$$\int_{-1}^{+1}\left[\frac{(\xi^2+1)(1-\eta^2)}{(\xi^2-\eta^2+1)}\right]^{m/2} y_m[c(\xi^2-\eta^2+1)^{1/2}]S_{ml}^{(1)}(-ic,\eta)d\eta,\quad l-m\ even. \tag{16}$$

$$R_{ml}^{(2)}(-ic,i\xi)=\frac{(-1)^{(l-m-1)/2}(2m+3)}{2^{m+1}m!d_1(-ic\,|\,ml)}\times$$
$$\int_{-1}^{+1}\frac{[(\xi^2+1)(1-\eta^2)]^{m/2}}{(\xi^2-\eta^2+1)^{(m+1)/2}}\xi\eta\, y_{m+1}[c(\xi^2-\eta^2+1)^{1/2}]S_{ml}^{(1)}(-ic,\eta)d\eta,\quad l-m\ odd. \tag{17}$$

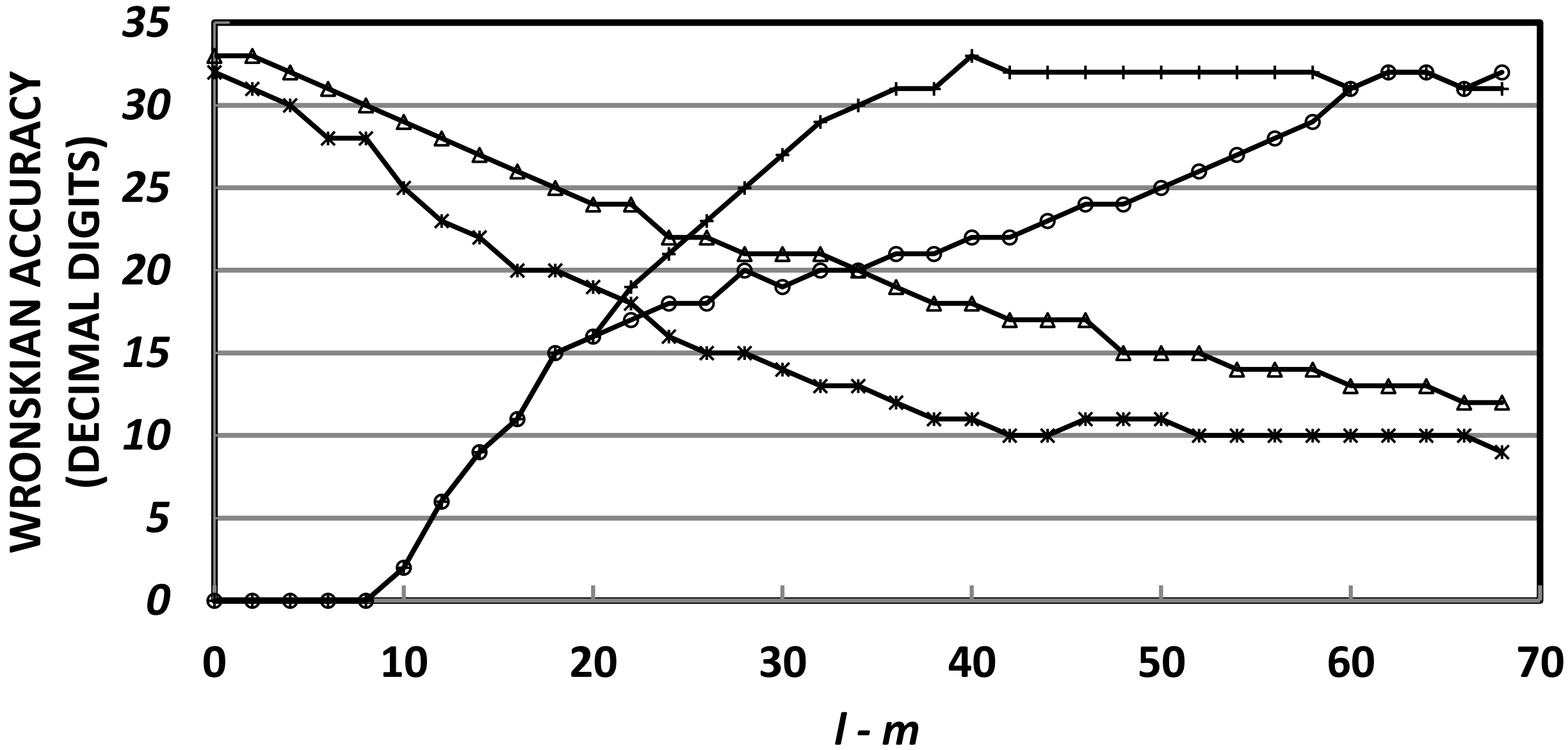

Fig. 8: Wronskian accuracy for two parameter sets ($\xi$, $c$, $m$) using the traditional Legendre expression: + (0.1, 60, 0); o (0.5, 30, 10) and using the Baber and Hasse expression: * (0.1, 60, 0); Δ (0.5, 60, 0)..

It is convenient to define $z=c(\xi^2-\eta^2+1)^{1/2}$ and a window function $F_m(\xi,\eta)$ given by

$$F_m(\xi,\eta)=\left[\frac{(\xi^2+1)(1-\eta^2)}{\xi^2-\eta^2+1}\right]^{m/2}. \tag{18}$$

Expanding $S_{ml}^{(1)}$ in (16) and (17) in terms of associated Legendre functions results in:

$$R_{ml}^{(2)}(-ic,i\xi)=B_{ml}^{(a)}(c){\sum_{n=0}^{\infty}}' d_n(-ic\,|\,ml)\,I_{mn}^{(a)}(c,\xi),\quad l-m\ even, \tag{19}$$

$$R_{ml}^{(2)}(-ic,i\xi)=\xi B_{ml}^{(b)}(c)\sum_{n=1}^{\infty}{}' d_n(-ic\,|\,ml)\,I_{mn}^{(b)}(c,\xi),\quad l-m\ odd, \tag{20}$$

where $B_{ml}^{(a)}(c)$ is the leading coefficient in (16), $B_{ml}^{(b)}(c)$ is the leading coefficient in (17), and

$$I_{mn}^{(a)}(c,\xi)=\int_{-1}^{+1}F_m(\xi,\eta)\,y_m(z)P_{m+n}^m(\eta)d\eta,\quad l-m\ even, \tag{21}$$

$$I_{mn}^{(b)}(c,\xi)=c\int_{-1}^{+1}\left[F_m(\xi,\eta)/z\right]\eta\, y_{m+1}(z)P_{m+n}^m(\eta)d\eta,\quad l-m\ odd. \tag{22}$$

One obtains corresponding expressions for the first derivatives of $R_{ml}^{(2)}$ with respect to $\xi$ from (16) and (17) by differentiating, utilizing standard recursion relations for the spherical Neumann functions, and collecting terms. This gives:

$$\begin{aligned}\frac{dR_{ml}^{(2)}}{d\xi}(-ic,i\xi)&=\frac{m\xi}{\xi^2+1}R_{ml}^{(2)}(-ic,i\xi)-\frac{(-1)^{(l-m)/2}(2m+1)c\xi}{2^{m+1}m!\,d_0(-ic\,|\,ml)}\times\\&\int_{-1}^{+1}\frac{\left[(\xi^2+1)(1-\eta^2)\right]^{m/2}}{(\xi^2-\eta^2+1)^{(m+1)/2}}\,y_{m+1}[c(\xi^2-\eta^2+1)^{1/2}]S_{ml}^{(1)}(-ic,\eta)d\eta,\quad l-m\ even,\end{aligned} \tag{23}$$

$$\begin{aligned}\frac{dR_{ml}^{(2)}}{d\xi}(-ic,i\xi)&=\frac{(m+1)\xi^2+1}{\xi(\xi^2+1)}R_{ml}^{(2)}(c,\xi)-\frac{(-1)^{(l-m-1)/2}(2m+3)c\xi^2}{2^{m+1}m!\,d_1(-ic\,|\,ml)}\times\\&\int_{-1}^{+1}\frac{[(\xi^2+1)(1-\eta^2)]^{m/2}}{(\xi^2-\eta^2+1)^{(m+2)/2}}\eta\, y_{m+2}[c(\xi^2-\eta^2+1)^{1/2}]S_{ml}^{(1)}(-ic,\eta)d\eta,\quad l-m\ odd.\end{aligned} \tag{24}$$

Replacing $S_{ml}^{(1)}$ with its expansion in (1) results in:

$$\frac{dR_{ml}^{(2)}}{d\xi}(-ic,i\xi)=\frac{m\xi}{\xi^2+1}R_{ml}^{(2)}(c,\xi)+c\xi B_{ml}^{(a)}\sum_{n=0}^{\infty}{}' d_n(-ic\,|\,ml)I_{mn}^{(c)}(c,\xi),\quad l-m\ even, \tag{25}$$

$$\frac{dR_{ml}^{(2)}}{d\xi}(-ic,i\xi)=\frac{(m+1)\xi^2+1}{\xi(\xi^2+1)}R_{ml}^{(2)}(-ic,i\xi)+c\xi^2B_{ml}^{(b)}\sum_{n=1}^{\infty}{}' d_n(-ic\,|\,ml)I_{mn}^{(d)}(c,\xi),\quad l-m\ odd, \tag{26}$$

where

$$I_{mn}^{(c)}(c,\xi)=\int_{-1}^{+1}\left[F_m(\xi,\eta)/z\right]y_{m+1}\left[z\right]P_{m+n}^m(\eta)d\eta,\quad l-m\ even, \tag{27}$$

$$I_{mn}^{(d)}(c,\xi)=\int_{-1}^{+1}\left[F_m(\xi,\eta)/z^2\right]\eta\, y_{m+2}\left[z\right]P_{m+n}^m(\eta)d\eta,\quad l-m\ odd. \tag{28}$$

The required integrals $I_{mn}(c,\xi)$ have an integrand that is symmetric about $\eta$ = 0. They can be computed using Gauss quadrature over positive values of $\eta$ and doubling the result. However, one must be careful to increase the density of quadrature points near $\eta$ = 1 when $\xi$ approaches zero because of the singularity of the spherical Neumann functions at $z$ = 0. The integrals tend to decrease in magnitude as $l$ increases, the decrease accompanied by loss of accuracy from increasing subtraction error. This causes a decrease in accuracy in both $R_{ml}^{(2)}$ and its first derivative as $l$ - $m$ increases, although the decrease can be very gradual in many cases.

Relationships between the different integrals can be obtained through use of recursion relations for the associated Legendre functions. For example, replacing $\eta P_{m+n}^m(\eta)$ in the rhs of (28) with its equivalent in terms of $P_{m+n-1}^m(\eta)$ and $P_{m+n+1}^m(\eta)$ results in:

$$(2n+2m+1)I_{mn}^{(b)} = (n+2m)I_{m,n-1}^{(c)} + (n+1)I_{m,n+1}^{(c)}. \tag{29}$$

Thus the integrals $I_{mn}^{(b)}$ can be calculated directly from (29) instead of computing them using Gauss quadrature. Other derived relations are not as useful for calculating radial functions with a given order *m* since they relate integrals of one kind and order *m* to integrals of a second kind and order $m \pm 1$. These could, however, be useful when one is computing the radial functions for a range of *m* values.

Use of the integral expressions provides accurate function values over a wide range of parameters. It is especially useful when $\xi \leq 0.5$, c is large, and *l* - *m* is less than the breakpoint. Here both of the Legendre expressions usually suffer too much subtraction error to allow accurate results. Figure 9 shows the Wronskian accuracy when using the integral expressions for several parameter sets. 128-bit arithmetic is used. Nearly full accuracy is obtained over a wide range of values of *l* - *m* when *m* is small, even when *c* is very large. When *m* is not small, the accuracy is reduced for lower values of *l* - *m*, the reduction increasing as *c* increases. The lowest accuracy occurs at *l* = *m*. However, even for values of *c* as large as 400 one obtains at least 10 decimal digits of accuracy. When *m* is much larger that the values shown here, accuracy at lower values of *l* - *m* is somewhat larger than shown here. For the case of $\xi = 0.1$ and *c* = 400, the accuracy at *l* = *m* is 13 digits for *m* = 200 and 18 digits for *m* = 250. When $\xi = 0.5$ and *c* = 400, the accuracy at *l* = *m* is 16 digits for m = 200 and 20 digits for m = 250. The falloff in accuracy at higher values of *l* - *m* seen in fig. 9 is unimportant since the traditional Legendre expression provides essentially full accuracy there.

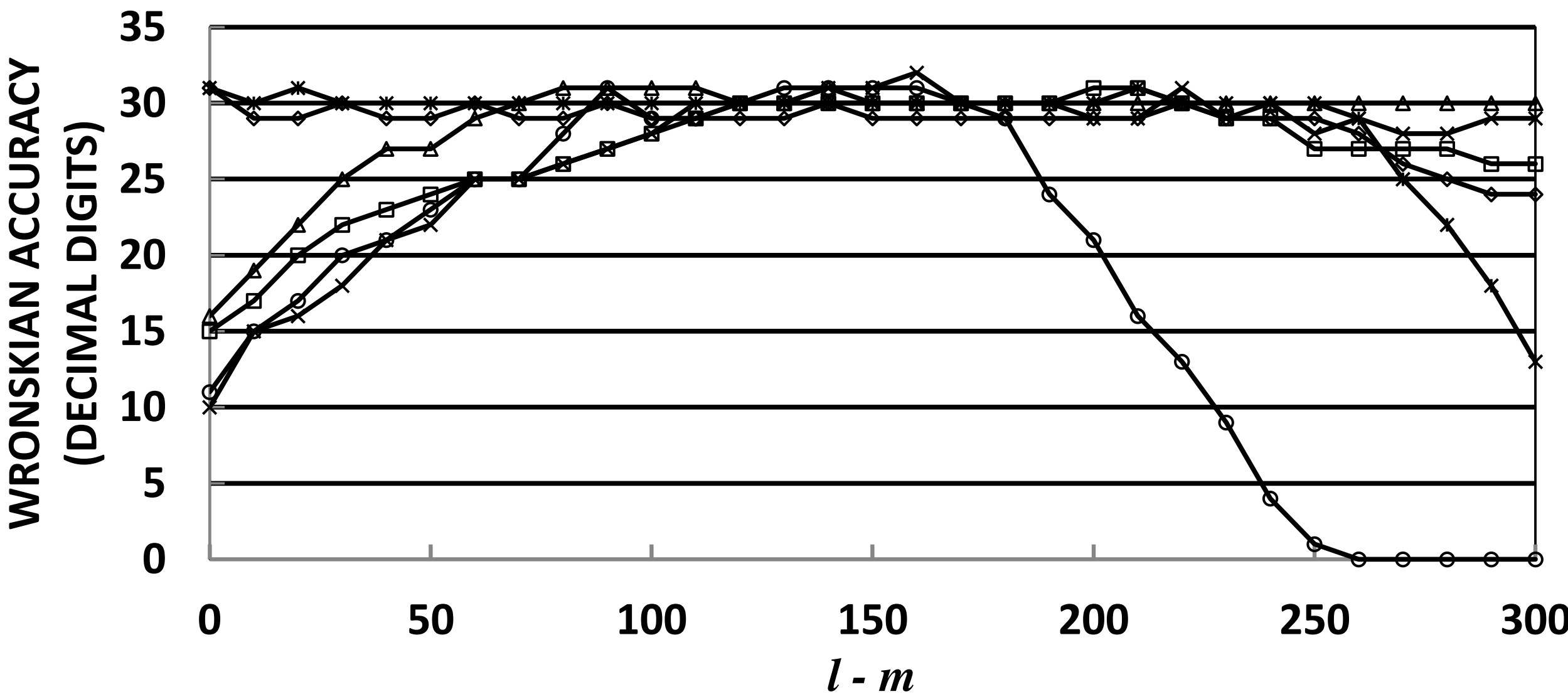

Fig. 9: Wronskian accuracy in decimal digits obtained using the integral method to calculate $R_{ml}^{(2)}(-ic, i\xi)$ for selected parameters ($\xi$, *c, m*): ◊ (0.1, 400, 0); □ (0.1, 400, 50); Δ (0.1, 300, 100); × (0.1, 400, 100); * (0.5, 400, 0); o (0.5, 400, 100).

## 11 Calculation of $R_{ml}^{(2)}(-ic,i\xi)$ using Neumann function expressions

Consider first the traditional Neumann function expression given in (8). This expression has the advantage that the denominator term is the corresponding angular function evaluated at $\eta$ = 1, which is numerically robust with no subtraction error. However, the numerator term is asymptotic and not absolutely convergent for any finite value of $c\xi$. It can often provide accurate values, especially when $\xi$ is not small. To evaluate (8) one takes the partial sum of the series up to and including the term where the magnitude of the relative contribution is smaller than $10^{-ndec}$. Here *ndec* is the number of decimal digits available in the arithmetic used in the calculations. The corresponding expression for the derivative of the radial function behaves similarly to (8). Sometimes the relative contribution never gets as small as $10^{-ndec}$. In that case the series is truncated at the term where the relative contribution is minimum. Figure 10 shows the accuracy obtained for $R_{ml}^{(2)}$ and its first derivative calculated in 128-bit arithmetic using (8) for selected parameter sets ($\xi$, *c, m*). The loss in accuracy with increasing *l* - *m* seen here is due to growing subtraction error. For all of these examples, the last term taken in the series had a relative contribution less than $10^{-33}$. The rate of loss in accuracy with increasing *l* - *m* grows as *c* decreases, as $\xi$ decreases and as *m* increases.

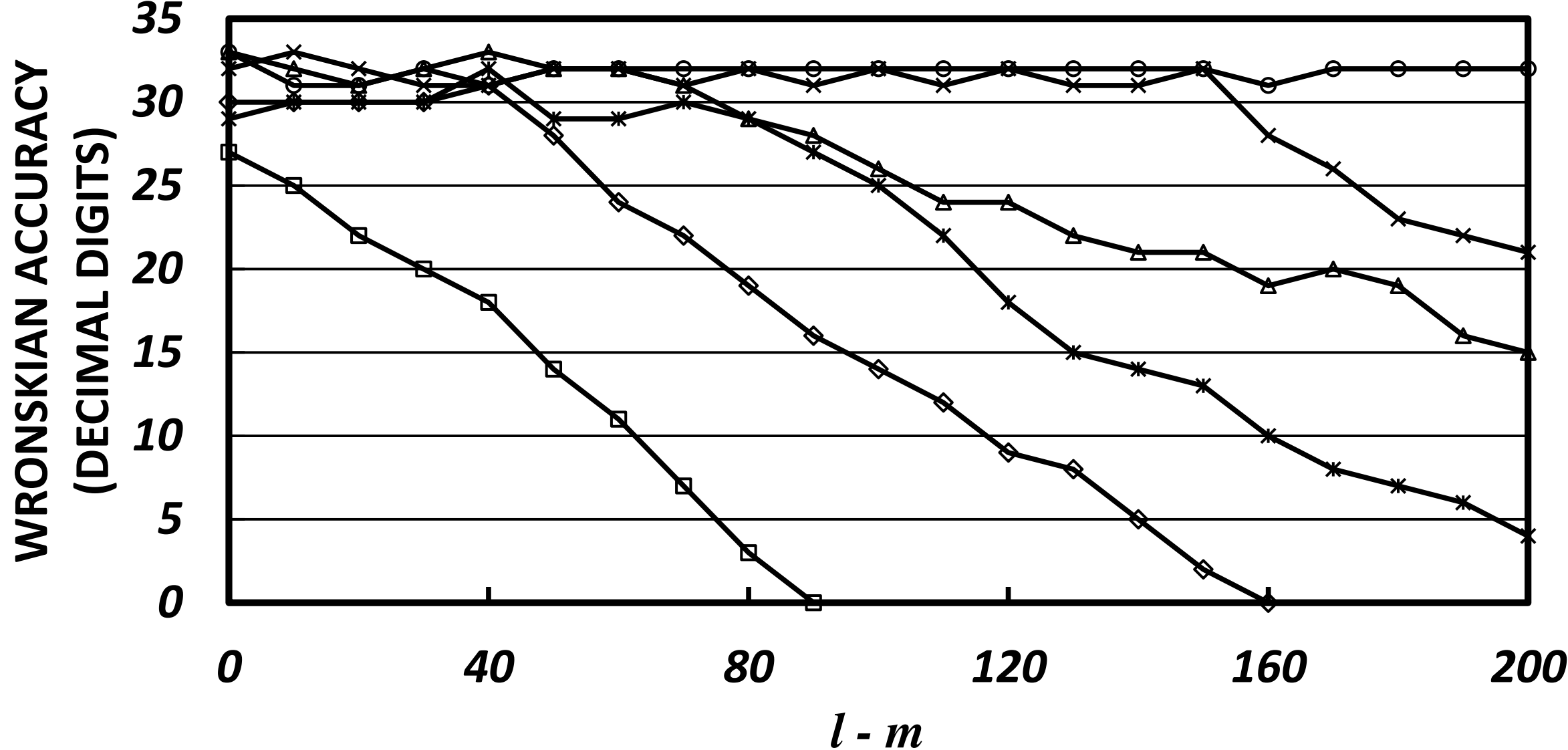

Fig. 10: Wronskian accuracy in decimal digits obtained using the traditional Neumann function expression to calculate $R_{ml}^{(2)}(-ic,i\xi)$ for selected parameters ($\xi$, *c, m*): ◊ (1.1, 40, 0); □ (1.1, 100, 40); Δ (1.5, 40, 0); × (1.5, 100, 0); * (1.5, 100, 50); o (1.5, 200, 0).

The number of terms required to achieve the accuracies shown here increases significantly with increasing *l* - *m*. For example, when $\xi$ = 1.1, *c* = 100, and m = 0, 77 terms are used at *l* = 0 and 232 terms are required at *l* = 30. The traditional Neumann expression becomes less useful when $\xi$ is less than unity. Here, a very large number of terms are often required and the relative contribution of successive terms may not decrease to the desired minimum. This is especially true when *m* is large. For example, when $\xi$ = 0.6, *c* = 300 and *m* =50, there is basically no minimum in relative contributions and the series fails to provide any accuracy.

The alternative Neumann function expressions (9) and (10) obtained when $\eta$ has been set equal to 0 are very useful. The numerator sums behave as if they were not asymptotic. They are well-behaved and converge to the desired accuracy, even at high values of $c$ and low values of $l$ - $m$. There is no evidence of the series beginning to diverge as further terms are added, even when tens of thousands of additional terms are taken in the series. This is true for values of $\xi$ as low as 0.01, although the number of terms required is much larger for lower values of $\xi$. Now the denominator in (9) and (10) does suffer the subtraction error of the corresponding angular function at $\eta = 0$. And the numerator sums suffer essentially the same subtraction error. As seen in figure 3, the error is greatest at $l = m$ and decreases to zero with increasing $l$ - $m$. For given $l$ - $m$ it increases with increasing $c$ and decreases with increasing $m$. Figure 2 shows this behavior for the case $l = m$. The subtraction error restricts the use of (9) and (10) to values of $l$ - $m$ that are large enough so that the numerator and denominator achieve the desired accuracy for the radial functions.

Figure 11 shows the accuracy of $R_{ml}^{(2)}$ and its first derivative obtained using (9) and (10)

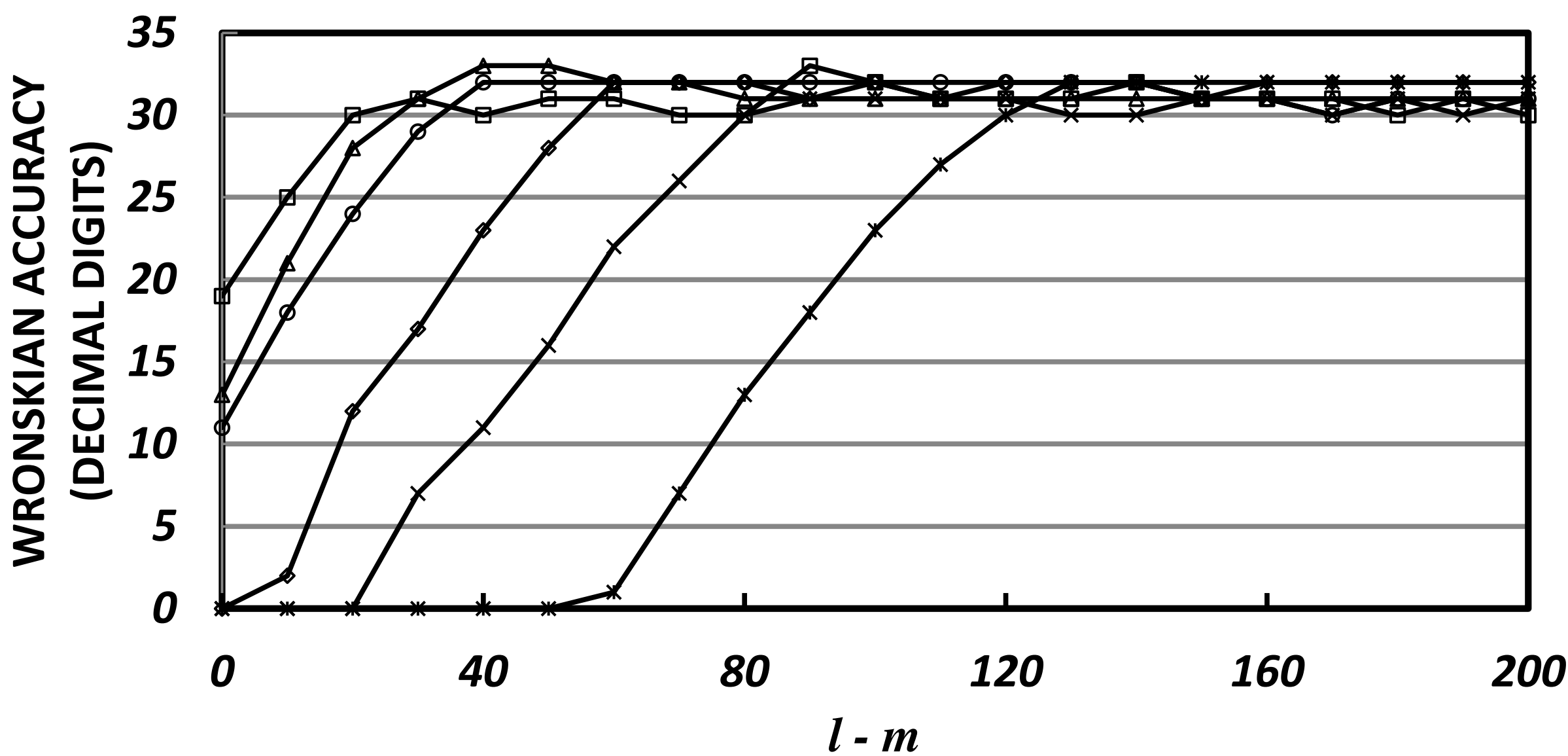

Fig. 11: Wronskian accuracy in decimal digits obtained using the Neumann function expression with $\eta = 0$ to calculate $R_{ml}^{(2)}(-ic, i\xi)$ for selected parameters ($\xi$, $c$, $m$): ◊ (0.5, 100, 0); □ (0.8, 100, 40); Δ (1.0, 50, 0); × (1.0, 200, 40); * (1.5, 200, 0); o (1.5, 200, 100).

for selected parameter sets ($\xi$, $c$, $m$). The accuracies shown here do not depend noticeably on the value of $\xi$ but only on the values for $c$ and $m$. They would remain virtually unchanged for any value of $\xi \geq 0.01$. For $\xi \geq 1$, the combination of the traditional Neumann function expression and the Neumann function expression with $\eta = 0$ often provides accurate results for all values of $l$ - $m$.

When $c$ is very large and $l$ - $m$ is below the breakpoint, sometimes the methods described above are unable to provide sufficiently accurate values for the radial functions of the second kind and their first derivatives. This is especially true when $m$ is neither small nor extremely large. The use of a variable $\eta$ method similar to that described above in section 6 can often help here when $\xi$ is greater than about 0.05. It can bridge the gap in $l$ - $m$ where eigenvalue pairing is no longer sufficient and the traditional Neumann expression (8) (or the integral method) fails to provide the desired accuracy and the expression with $\eta = 0$ begins to do so. Figure 11 shows accuracies obtained for several values of $\eta$ when $\xi = 0.6$, $c = 300$ and $m = 50$. The results shown for $\eta = 0$ are those obtained using (9) and (10). Just four different values of $\eta$ provide a minimum of 20 decimal digits of Wronskian accuracy for all $l - m$. Including the results for $\eta = 0.8$ increases the minimum accuracy to 23 digits. Of course one does not know in advance how well particular $\eta$ values work. More values of $\eta$ will usually be tried in order to obtain the desired accuracy. Although the Neumann function series are asymptotic, the minimum contribution term is often less than $10^{-33}$.

The variable $\eta$ method provides even better results for larger values of $c$. For $\xi \geq 0.5$ and $c$ near or well above 1000, it can provide nearly full accuracy for lower values of $l - m$ where other methods fail. And it usually does so with only a few values of $\eta$. As an example, consider the case where $\xi = 1.0$, $c = 1000$, $m = 100$ and 28 decimal digits of accuracy are requested. As seen in Fig. 12, the single value of $\eta = 0.954$ provided a minimum of 28 digits of accuracy from $l = 100$ to $l = 624$, where the $\eta = 0$ expressions began providing at least 30 digits of accuracy. The addition of the variable $\eta$ method allows one to calculate radial function values with more than 16 digits of accuracy for values of c ranging up to 5000 or more when $\xi \geq 0.5$.

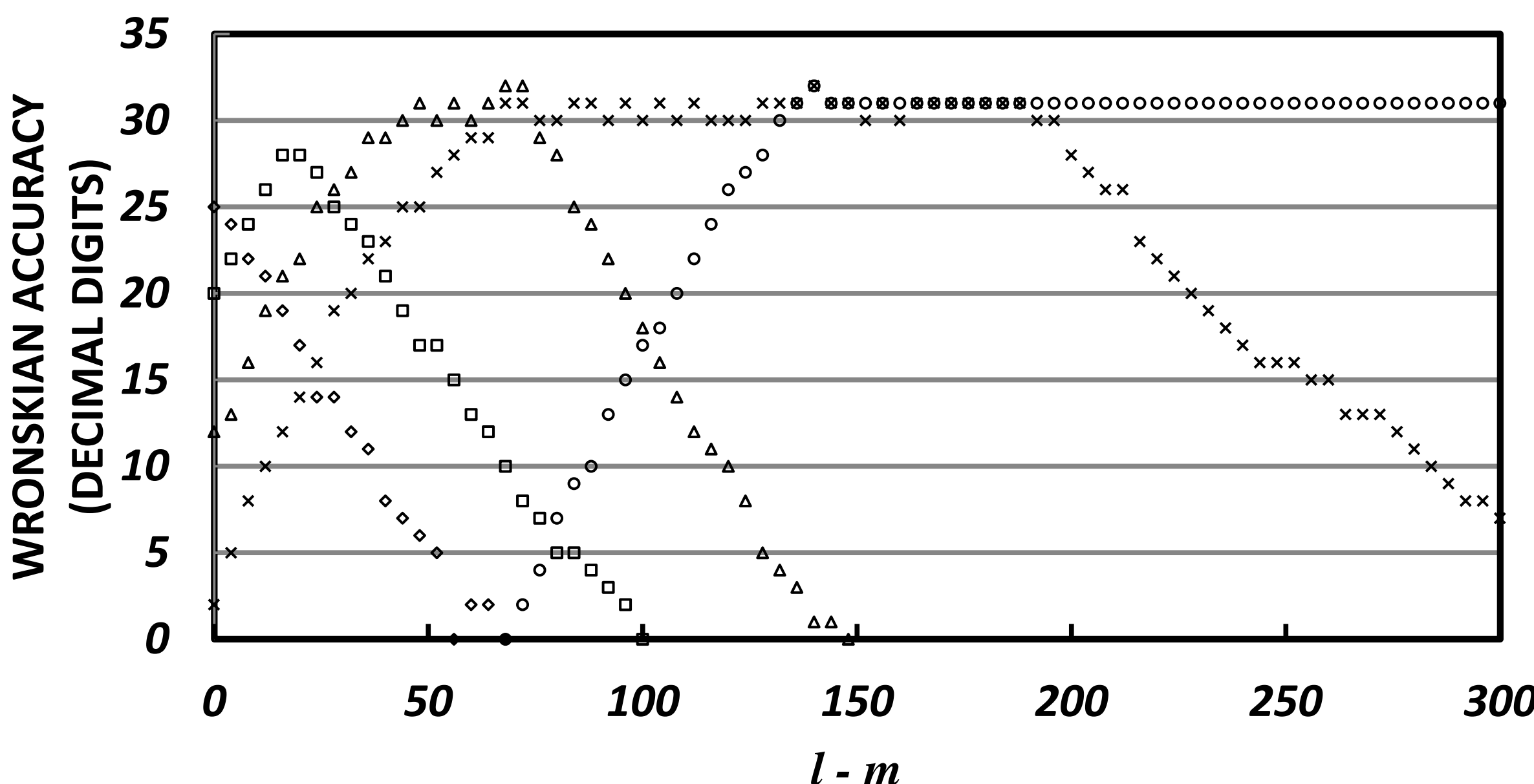

Fig. 12: Wronskian accuracy in decimal digits obtained for $\xi = 0.6$, $c = 300$, and $m = 50$ using the variable $\eta$ method with selected values of $\eta$: ◊ (0.8); □ (0.7); Δ (0.6); × (0.5); o (0.0).

The variable $\eta$ method is implemented as follows. At the lowest value of $l - m$ where the desired accuracy is not achieved using other methods, the value of $\eta$ is decreased in steps from unity and the radial functions of the second kind are calculated at each step. The accuracy initially tends to increase with decreasing $\eta$ as the numerator series becomes more accurate, although it may take several steps before any significant increase is obtained. Ideally, the desired accuracy is achieved after one of the steps. The associated value of $\eta$ for that step is then used for the next value of $l - m$. It continues to be used for progressively higher values of $l - m$ until the accuracy again falls below the desired minimum. Then $\eta$ is again decreased in steps until the desired accuracy is achieved for that value of $l - m$. Typically only one or two steps are needed here. The process is repeated until the $\eta = 0$ expression offers sufficient accuracy. Sometimes no value of $\eta$ provides the desired accuracy, even when the process is continued until the denominator becomes less accurate than the numerator series and the accuracy starts to decrease. The best $\eta$ for this value of $l$ - $m$ is the one used for the previous step, which is then used for the next value of $l$ - $m$. The process is continued until the $\eta = 0$ expressions offers the desired accuracy and is used instead. When $\xi \leq 0.99$, the traditional Legendre expressions will replace the $\eta = 0$ expressions when $l$ - $m$ is large enough so that they provide the desired accuracy. As in section 6 it is convenient to use steps in $\theta = \arccos(\eta)$. A step size of about 0.1 radian is used when $\xi > 0.4$ and about 0.05 radian for lower values of $\xi$.

## 12 Calculation of radial functions for $\xi = 0$

When $\xi = 0$, accurate values for $R_{ml}^{(1)}$ are given when $l$ - $m$ is even by the $d_0$ term in (8) since the remaining terms vanish. $R_{ml}^{(1)}$ for $l$ - $m$ odd is equal to zero since all of the terms vanish. The first derivative of $R_{ml}^{(1)}$ for $l$ - $m$ odd is given by the $d_1$ term in the derivative of (8) since the remaining terms vanish. The first derivative of $R_{ml}^{(1)}$ is equal to zero for $l$ - $m$ even since all of the terms vanish. The Wronskian can be used to obtain accurate values for $R_{ml}^{(2)}$ for $l$ - $m$ odd and for its first derivative when $l$ - $m$ is even. Here

$$
\begin{aligned}
R_{ml}^{(2)}(-ic,0) &= -\frac{1}{c\,\dfrac{dR_{ml}^{(1)}(-ic,0)}{d\xi}}, \qquad l-\mathrm{m}\ \ odd, \\
\frac{dR_{ml}^{(2)}(-ic,0)}{d\xi} &= \frac{1}{cR_{ml}^{(1)}(-ic,0)}, \qquad l-m\ \ even.
\end{aligned}
\tag{30}
$$

Values for $R_{ml}^{(2)}$ when $l$ - $m$ is even are obtained from the limiting form of (12) while values for its first derivative for $l$ - $m$ odd are given by the limiting form of the first derivative of (12). These limiting forms suffer the same subtraction errors as those of the angular function of the first kind when $\eta = 0$ as shown in Fig. 2. The resulting function values are progressively smaller relative to the corresponding non-zero function values for $R_{ml}^{(1)}$ or its first derivative. Thus their contribution to the solution of problems involving oblate spheroidal geometry is proportionately reduced along with their reduced accuracy. When the subtraction error is within 3 of the number of decimal digits available, the function values are set equal to zero.

## 13 A Fortran program to compute oblate spheroidal functions

A Fortran computer program called oblfcn [12] has been developed to calculate the oblate spheroidal functions. Oblfcn is available as either a stand-alone program or as a subroutine. It performs calculations in either double precision (64 bit) arithmetic or in quadruple precision (128 bit) arithmetic. The choice of arithmetic is controlled by a module called param that is compiled prior to compiling oblfcn and that sets the kind parameter equal to either 8 for real*8 or to 16 for real*16. Calculation options include (1) radial functions of the first kind and their first derivatives, (2) radial functions of both the first and second kind and their first derivatives, (3) angular functions of the first kind, and (4) angular functions of the first kind and their first derivatives. If desired both radial and angular functions can be calculated during the same run.

Oblfcn provides an estimate of the number of accurate digits in the angular functions and their first derivatives based on the subtraction errors involved in their calculation. It also provides an estimate of the number of accurate digits in $R_{ml}^{(2)}$ and its first derivative. Both $R_{ml}^{(1)}$ and its first derivative are always highly accurate. The estimate of accuracy is often based on the Wronskian. When (11) is used, it is based on the number of digits of agreement between low-order paired eigenvalues. When $\xi = 0$ , the estimate of accuracy is based on the subtraction error involved in the calculation of either $R_{ml}^{(2)}$ when $l$ - $m$ is even or the first derivative of $R_{ml}^{(2)}$ when $l$ - $m$ is odd. The output of oblfcn includes diagnostic files (one for radial functions and one for angular functions) including information such as the number of terms both available and used in the various series.

Oblfcn calculates angular functions of the first kind using (1) and gives them unit norm or normalizes them using (5), depending on the input parameter iopnorm. It does so for a specified number of $l$ values beginning with $l = m$. In the stand-alone program it allows either $\eta$ or $\theta = \arccos(\eta)$ arguments and computes angular functions for a range of arguments determined by a first value, an increment and the number of arguments desired. It does so for a range of $m$ values determined by a first value, an increment and the number of $m$ values desired. The resulting function values are given as a characteristic with a magnitude between 1.0 and 10.0 and an integer exponent iexp that denotes the power of 10 for the factor 10.0**(iexp). In the subroutine version angular functions values are obtained for a single value of $m$ and an input vector of $\eta$ values. The resulting values are not provided as a characteristic and exponent but as a single number. It is expected that the user will choose unit norm to avoid any potential overflow problems for very high values of $m$.

Oblfcn calculates radial function for a single input value of $\xi$ and for a specified number of $l$ values beginning with $l = m$. In the stand-alone program it does so for a range of $m$ values while in the subroutine version it does so only for a single value of $m$. Oblfcn calculates the radial functions of the first kind using the traditional Bessel function expression together with the variable $\eta$ method described above. The radial functions of the first kind are nearly fully accurate (unless near a root). Oblfcn calculates radial functions of the second kind using either the pairing of low-order eigenvalues and (11), traditional $\eta = 1$ Neumann function expression (8), the alternative $\eta = 0$ Neumann function expressions (9) - (10), the traditional associated Legendre function expression (12), the Baber and Hasse Legendre function expression (14), the variable $\eta$ method (7), or the integral expressions in (16) and (17). The methods used in oblfcn are based on the input parameters and the integer minacc that specifies the number of accurate decimal digits that are desired. Minacc is set equal to 10 for 64-bit arithmetic and to 15 for 128-bit arithmetic. The value for 128-bit arithmetic can be changed if desired, especially if higher accuracy is desired for input parameters where it can be achieved. It is advised to leave minacc set to 10 for 64-bit arithmetic.

The methods used in Oblfcn to calculate the radial functions of the second kind are based on the input parameters and the desired minimum accuracy. Oblfcn starts at $l = m$ with the use of paired eigenvalues if the pairing is sufficient to provide the desired accuracy and continues with increasing $l$ - $m$ until the pairing is insufficient. For values of $c \leq 20$ for real*8 or $c \leq 60$ for real*16, oblfcn then uses the alternative $\eta = 0$ Neumann function expressions for $\xi > 0.99$ and the traditional Legendre function expression for $\xi \leq 0.99$. For intermediate values of $c$, it may also use the Baber and Hasse Legendre expression at lower values of $l$ - $m$ where the traditional Legendre function expression is insufficiently accurate.

For larger values of c and smaller values of $\xi$, oblfcn switches to the integral method after the paired eigenvalue method and continues until the traditional Legendre function expansion provides sufficient accuracy and is used for all higher values of $l$ - $m$ that are desired . For larger values of both $c$ and $\xi$, it switches to one of the Neumann function expansions and continues until the $\eta = 0$ Neumann expansion provides the desired accuracy and is used for all higher values of $l$ - $m$ that are desired.

The various series involved in oblfcn are computed in a way to avoid potential overflow and underflow in the calculations. First the expansion coefficients and the expansion functions are calculated as ratios using appropriate recursion relations. See comment statements in the appropriate subroutine in oblfcn where they are calculated. Bessel and Neumann function series are summed starting with the $n = l - m$ term while Legendre function and integral method series are summed starting with the lowest term. The function and coefficient values for the first term are factored out of the expansion and the first term is set equal to unity. Summation is performed using ratios to obtain the next terms in the series For Bessel and Neumann functions expansions summations are taken both forward and backward. The resulting value is then multiplied by the first term to obtain the desired sum. Here the relevant Bessel, Neumann, or Legendre function for the first term taken has been computed previously by forward multiplication of ratios starting with known values for the lowest two values. For example, this would be $j_0$ and $j_1$ for Bessel functions. During forward multiplication, power 10 exponents are stripped out of the product at each step to avoid either underflow or underflow. This results in a characteristic with magnitude between 1.0 and 10.0 and an exponent denoting the corresponding power of 10. Radial function values are stored as both a characteristic and an exponent. This allows oblfcn to provide results at high values of $l$ - $m$ where the radial functions of the first kind would otherwise underflow and the radial functions of the second kind would overflow.

Oblfcn was tested extensively using 64-bit arithmetic for values of $m$ up to 1000, for $c$ up to 10000 and for a wide range of values of $\xi$ down to and including zero. Accuracies of 8 or more decimal digits were usually obtained. Accuracies less than 8 digits can occur when $\xi$ is near a root. They can also occur when $c$ is very large, $\xi$ is near or equal to zero and $l$ - $m$ is small. In this case either $R_{ml}^{(2)}$ or its first derivative has a reduced magnitude along with its reduced accuracy and is not likely to have an adverse impact on the solution to problems using these function values. Use of 128-bit arithmetic extends the accuracy to 15+ digits. It is noted that computation times using 128-bit arithmetic are up to 50 or more times larger than those using 64-bit arithmetic. It is likely that the accuracy obtained using 64-bit arithmetic is sufficient for most applications.

## 14 Summary

A new procedures is presented which provides values for the oblate radial function of the first kind and its first derivative that are nearly fully accurate for virtually all parameter ranges. Several alternative expressions combined with traditional ones for calculating the oblate spheroidal radial functions of the second kind and their first derivatives provide accurate values over extremely wide parameter ranges even when using 64-bit arithmetic. This paper also briefly describes a Fortran computer program oblfcn that incorporates these new procedures. Oblfcn can be run in either 64-bit or 128-bit arithmetic A listing of oblfcn together with sample output is freely available in text format on the web site listed in [12]. Both a stand-alone program and a subroutine version are provided.